\def\yes{\if00}
\def\no{\if01}
\def\iftwelvept{\yes}
\def\ifusepdf{\no}
\def\ifusepsfont{\no}

\iftwelvept
\documentclass[leqno,12pt]{amsart}
\else
\documentclass[leqno]{amsart}
\fi
\usepackage{amsfonts}
\usepackage{amssymb}
\usepackage{amscd}
\usepackage{latexsym}
\usepackage{verbatim}
\ifusepdf
\usepackage{hyperref}
\else\fi
\ifusepsfont
\usepackage[T1]{fontenc}
\usepackage{times}
\else\fi

\iftwelvept
\else\fi

\newtheorem{Theorem}{Theorem}[subsection]
\newtheorem{Proposition}[Theorem]{Proposition}
\newtheorem{Lemma}[Theorem]{Lemma}
\newtheorem{Corollary}[Theorem]{Corollary}
\newtheorem{Claim}{Claim}[Theorem]

\theoremstyle{definition}
\newtheorem{Definition}[Theorem]{Definition}
\newtheorem{Remark}[Theorem]{Remark}
\newtheorem{Example}[Theorem]{Example}

\newtheorem{Question}[Theorem]{Question}

\renewcommand{\theTheorem}{\arabic{section}.\arabic{subsection}.\arabic{Theorem}}
\renewcommand{\theClaim}{\arabic{section}.\arabic{subsection}.\arabic{Theorem}.\arabic{Claim}}
\renewcommand{\theequation}{\arabic{section}.\arabic{subsection}.\arabic{Theorem}.\arabic{Claim}}

\newcommand{\ZZ}{{\mathbb{Z}}}
\newcommand{\QQ}{{\mathbb{Q}}}
\newcommand{\RR}{{\mathbb{R}}}
\newcommand{\CC}{{\mathbb{C}}}
\newcommand{\PP}{{\mathbb{P}}}
\newcommand{\OO}{{\mathcal{O}}}
\newcommand{\Aff}{{\mathbb{A}}}

\newcommand{\Aut}{\operatorname{Aut}}
\newcommand{\Spec}{\operatorname{Spec}}
\newcommand{\Pic}{\operatorname{Pic}}
\newcommand{\ord}{\operatorname{ord}}
\newcommand{\Div}{\operatorname{Div}}
\newcommand{\Supp}{\operatorname{Supp}}
\newcommand{\zero}{\operatorname{div}}
\newcommand{\Rat}{\operatorname{Rat}}

\newcommand{\calF}{{\mathcal{F}}}
\newcommand{\calX}{{\mathcal{X}}}
\newcommand{\calY}{{\mathcal{Y}}}
\newcommand{\calZ}{{\mathcal{Z}}}
\newcommand{\calW}{{\mathcal{W}}}
\newcommand{\calL}{{\mathcal{L}}}
\newcommand{\calM}{{\mathcal{M}}}
\newcommand{\frakP}{{\mathfrak{P}}}
\newcommand{\id}{\operatorname{id}}
\newcommand{\Image}{\operatorname{Im}}

\newcommand{\chern}{\operatorname{c}}
\newcommand{\achern}{\operatorname{\widehat{c}}}
\newcommand{\adeg}{\operatorname{\widehat{deg}}}
\newcommand{\topo}{\operatorname{top}}

\newcommand{\Proof}{{\sl Proof.}\quad}

\newcommand{\QED}{{\unskip\nobreak\hfil\penalty50\quad\null\nobreak\hfil
{$\Box$}\parfillskip0pt\finalhyphendemerits0\par\medskip}}
\newcommand{\rest}[2]{\left.{#1}\right\vert_{{#2}}}

\begin{document}

\title[canonical heights]%
{Canonical heights, invariant currents, and dynamical systems
of morphisms associated with line bundles}
\author{Shu Kawaguchi}
\address{Department of Mathematics, Faculty of Science,
Kyoto University, Kyoto, 606-8502, Japan}
\email{kawaguch@math.kyoto-u.ac.jp}
\begin{abstract}
We construct canonical heights of subvarieties 
for dynamical systems of several morphisms associated with line bundles 
defined over a number field, and study some of their properties. 
We also construct invariant currents for such systems over $\CC$. 
\end{abstract} 

\maketitle

\section*{Introduction}
\renewcommand{\theTheorem}{\Alph{Theorem}} 
Let $X$ be a projective variety over a field $K$ and 
$f_i: X \to X$ ($i=1, \cdots, k$) morphisms over $K$.  
Let $L$ be a line bundle on $X$, and $d > k$ a real number.  
We say that a pair $(X; f_1, \cdots, f_k)$ is a 
{\em dynamical system of $k$ morphisms  
over $K$ associated with $L$ of degree $d$} 
if $\bigotimes_{i=1}^{k} f_i^{*}(L) \simeq L^{\otimes d}$ holds.  
The purpose of this paper to construct 
canonical heights of subvarieties for such systems 
when $K$ is a number field 
and study some of their properties, 
and construct 
invariant currents for such systems when $K$ is $\CC$. 
We also remark on distribution of points of small heights 
for Latt\`es examples, which are certain endomorphisms of $\PP^N$.  

Firstly, we explain canonical heights. 
Weil heights play one of the key roles in Diophantine geometry, and 
particular Weil heights that enjoy nice 
properties (called ``canonical'' heights) are 
sometimes of great use. 

Over abelian varieties $A$ defined over a number field $K$, 
N\'eron and Tate constructed height functions 
(called canonical height functions or N\'eron--Tate height functions) 
$\widehat{h}_L: A(\overline{K}) \to \RR$ with respect 
to symmetric ample line bundles $L$ which enjoys nice properties.  
More generally, in \cite{CS} Call and Silverman constructed 
canonical height functions on 
projective varieties $X$ defined over a number field 
which admit a morphism $f: X \to X$ with 
$f^{*}(L) \simeq L^{\otimes d}$ for some 
line bundle $L$ and some $d >1$. 
In another direction, 
Silverman \cite{SiK3} constructed canonical
height functions on certain K3 surfaces $S$ 
with two involutions $\sigma_1, \sigma_2$ 
(called Wheler's K3 surfaces, cf. \cite{Maz}) 
and developed an arithmetic theory
analogous to the arithmetic theory on abelian varieties. 

Regarding canonical heights of subvarieties of projective varieties, 
Philippon \cite{Ph}, Gubler \cite{G} and Kramer \cite{G} 
constructed canonical heights of subvarieties of 
abelian varieties. 
In \cite{Zh}, Zhang  constructed canonical 
heights of subvarieties of projective varieties $X$ 
which admit a morphism $f: X \to X$ with 
$f^{*}(L) \simeq L^{\otimes d}$ for some  
line bundle $L$ and some $d >1$. 
For Wheler's K3 surfaces, however, canonical heights 
of subvarieties seem not to have been constructed. 

Our first results are construction of canonical heights 
of subvarieties for dynamical system 
of morphisms associated with line bundles. 

\begin{Theorem}[cf. Theorem~\ref{thm:main},  
Proposition~\ref{prop:property:of:height}, 
Theorem~\ref{thm:canonical:height:of:subvarieties}, 
Theorem~\ref{thm:global:and:local:canonical:height}]
\label{thm:Intro}
Let $(X;$ $f_1, \cdots, f_k)$ be a dynamical system of 
$k$ morphisms over a number field $K$ 
associated with a line bundle $L$ of degree $d > k$. 
Then 
\begin{enumerate}
\item[(1)]
there exists a unique real-valued function 
\[
\widehat{h}_{L, \{f_1, \cdots, f_k\}}: X(\overline{K}) \to \RR
\]
with the following properties\textup{:}
\begin{enumerate}
\item[(i)]
$\widehat{h}_{L, \{f_1, \cdots, f_k\}}$ is a Weil height corresponding to $L$, 
i.e., $\widehat{h}_{L, \{f_1, \cdots, f_k\}} = h_L + O(1)$\textup{;} 
\item[(ii)]
$\sum_{i=1}^{k} \widehat{h}_{L, \{f_1, \cdots, f_k\}} (f_i(x)) = 
d \ \widehat{h}_{L, \{f_1, \cdots, f_k\}}(x)$ for all $x \in X(\overline{K})$. 
\end{enumerate}
\item[(2)]
Assume $L$ is ample. Then 
\begin{enumerate}
\item[(a)]
$\widehat{h}_{L, \{f_1, \cdots, f_k\}}(x) \geq 0$ 
for all $x \in X(\overline{K})$. 
\item[(b)]
Let $C(x) :=\{f_{i_1}\circ\cdots\circ f_{i_l}(x) 
\mid l \geq 0, 1\leq i_1, \cdots, i_l \leq k\}$ denote 
the forward orbit of $x$ under $\{f_1, \cdots, f_k\}$. 
Then $\widehat{h}_{L, \{f_1, \cdots, f_k\}}(x) = 0$ if and only if 
$C(x)$ is finite. 
\end{enumerate}
\item[(3)]
Assume that $L$ is ample, $X$ is normal and $f_1, \cdots, f_k$ 
are all surjective. Then one can define 
the canonical height 
$\widehat{h}_{L, \{f_1, \cdots, f_k\}}(Y) \geq 0$ for 
any subvariety $Y \subset X_{\overline{K}}$.  
\item[(4)]
Assume $X$ is normal. Then 
the function $\widehat{h}_{L, \{f_1, \cdots, f_k\}}$ decomposes into 
the sum of local canonical height functions. 
\end{enumerate}
\end{Theorem}

Consider Wheler's K3 surface $S$. 
It is shown that $(S;\sigma_1, \sigma_2)$ becomes 
a dynamical system of 
two morphisms associated with a certain ample line bundle 
of degree $4$. Applying Theorem~A to $(S; \sigma_1, \sigma_2)$, 
one obtains a height function 
$\widehat{h}_{L, \{\sigma_1, \sigma_2\}}: 
S(\overline{K}) \to \RR$. 
It will be shown that $(1+\sqrt{3})\widehat{h}_{L, \{\sigma_1, \sigma_2\}}$ 
coincides with the the canonical height function 
defined by Silverman. 
In this way, Theorem~A generalizes Call--Silverman and 
Zhang's construction of canonical heights 
for dynamical systems of one morphism 
to dynamical systems of $k$ morphisms, 
incorporating Silverman's canonical height 
functions on Wheler's K3 surfaces. 
Moreover, one has canonical heights of subvarieties 
of Wheler's K3 surfaces. 

Using Theorem~\ref{thm:Intro}, we can similarly 
construct canonical height functions (and canonical heights of subvarieties) 
on other K3 surfaces: K3 surfaces 
in $\PP^2\times\PP^2$ given by the intersection of two hypersurfaces 
of bidegrees $(1,2)$ and $(2,1)$; and K3 surfaces 
in $\PP^1\times\PP^1\times\PP^1$ given 
by hypersurfaces of tridegree $(2,2,2)$. 
For arithmetic properties of the last surfaces, 
see Mazur \cite{Maz} and Wang \cite{Wa}. 
Another example of a dynamical system is the projective space $\PP^N$ 
and several surjective morphisms $f_i: \PP^N \to \PP^N$ such 
that $\deg f_i \geq 2$ for some $i$. 
For their properties over $\CC$ in the case $N=1$ 
(dynamics of finitely generated rational semigroups),  
we refer to Hinkkanen and Martin \cite{HM}, 
which studies the theory of the dynamics 
of semigroups of rational functions on $\PP^1$. 
See \S\ref{subsec:examples} for other examples.

Secondly, we explain invariant currents. 
In \cite{HP}, Hubbard and Papadopol introduced 
an invariant current, called the Green current, 
for a morphism $f:\PP^N \to \PP^N$ of degree $d \geq 2$, 
which is a generalization of the Brodin measure (also called 
the equilibrium measure) of $N=1$. 
(More strongly, the Green current of $f$ is defined 
when $f:\PP^N \cdots\to \PP^N$ is an algebraically 
stable rational map, cf. \cite{Sib},~Th\'eor\`eme~1.6.1.) 
Properties of the Green current of $f$ have since been deeply 
studied by many authors, see for example the survey of Sibony 
\cite{Sib} and the references therein. 

Our second results are the existence of 
$(1,1)$-currents with nice properties, 
which we call invariant currents, 
for dynamical systems of morphisms over $\CC$  
associated with line bundles. 

\begin{Theorem}[cf. Theorem~\ref{thm:admissible:metrics}, 
Theorem~\ref{thm:invariant:currents}]
\label{thm:Intro:2}
Let $(X; f_1, \cdots, f_k)$ be a dynamical system of $k$ morphisms over
$\CC$ associated with a line bundle $L$ of degree $d > k$. We assume that
$X$ is smooth and $L$ is ample.
\begin{enumerate}
\item[(1)]
Let $\eta_0 \in A^{1,1}(X(\CC))$ be a closed $C^{\infty}$ $(1,1)$-form
whose cohomology class coincides with $\chern_1(L)$. We inductively define
$\eta_n \in A^{1,1}(X(\CC))$ by
\[
\eta_{n+1} = \frac{1}{d} (f_1^*\eta_n + \cdots + f_k^*\eta_n).
\]
Then, $[\eta_n]$ converges to a closed positive $(1,1)$-current $T$. 
Moreover, $T$ is independent of the choice of $\eta_0$, and satisfies 
$f_1^*T + \cdots + f_k^*T = dT$.
\item[(2)]
The current $T$ admits a locally continuous potential. 
Indeed, there is a unique continuous metric 
$\Vert\cdot\Vert_{\infty}$, 
called the admissible metric \textup{(}cf. \cite{Zh}\textup{)}, 
on $L$ with $\Vert\cdot\Vert_{\infty}^d =
\varphi^*(f_1^*\Vert\cdot\Vert_{\infty} \cdots
f_k^*\Vert\cdot\Vert_{\infty})$, 
and $T$ is given by $\chern_1(L, \Vert\cdot\Vert_{\infty})$. 
\end{enumerate}
\end{Theorem}

When $f:\PP^N \to \PP^N$ is a morphism of degree $d \geq 2$, 
the current $T$ for $(\PP^N;f)$ is nothing but 
the Green current of $f$. Consider Wheler's K3 surface 
$S$. In this case, $\sigma_2\circ\sigma_1$ (resp. $\sigma_1\circ\sigma_2$) 
has a strictly positive topological entropy, 
and by the results of Cantat \cite{Ca} 
there exists a closed positive $(1,1)$-current $T^{+}$ 
(resp. $T^{-}$), unique up to scale, such that 
$(\sigma_2\circ\sigma_1)^*(T^{+}) = (7+4\sqrt{3}) T^{+}$ 
(resp. $(\sigma_1\circ\sigma_2)^*(T^{-}) = (7+4\sqrt{3}) T^{-}$). 
It will be shown that the current $T$ for the dynamical system 
$(S; \sigma_1, \sigma_2)$ is equal to $T^{+} + T^{-}$ with a 
certain normalization of a scale.

Since Szpiro, Ullmo and Zhang \cite{SUZ} proved equidistribution 
of small points for abelian varieties, and 
since heights of subvarieties and currents both with 
nice properties are constructed for $(X; f_1, \cdots, f_k)$, 
it is tempting to pose a question for $(X; f_1, \cdots, f_k)$ 
how the Galois orbit of $x_n$ is distributed with respect to  
$\mu:= \frac{T \wedge \cdots \wedge T}{\chern_1(L)^{\dim X}}$ 
($\dim X$-times) whenever there exists 
a sequence $\{x_n\}_{n=1}^{\infty}$ of $X(\overline{K})$ 
of small points. Using the methods of \cite{SUZ} and \cite{CL}, 
we remark that this question is true for Latt\`es examples, 
which are endomorphisms of $\PP^N$ for which the Green currents 
are on some non-empty (analytic) open subset smooth 
and strictly positive. 

The organization of this paper is as follows. 
In \S1, we construct canonical heights, study some of 
their properties, and give some examples. 
In \S2, we construct canonical heights of 
subvarieties of $X$ using adelic intersection theory developed 
by Zhang. 
In \S3, we construct invariant currents on $X(\CC)$. 
Then, we remark on distribution of points of small height 
for Latt\`es's examples. 
In \S4, we construct canonical 
local heights and see that the canonical height in \S1 
decomposes into the sum of these canonical local heights. 
The case $k=1$ is studied by Call and Silverman. 
In the appendix, we show finiteness of $\{f_1, \cdots, f_k\}$-periodic 
points of bounded degree in a more general setting. 

The author expresses his sincere gratitude to 
Professors Atsushi Moriwaki and Vincent Maillot 
for helpful conversations, to Professor Tetsuo Ueda for 
explaining Latt\`es examples of $\PP^N$. 

\medskip
\renewcommand{\theTheorem}{\arabic{section}.\arabic{subsection}.\arabic{Theorem}}
\section{Canonical heights}
\label{sec:canonical:heights}
\subsection{Quick review of heights}
\label{subsec:quick:review:of:heights}
In this subsection, we briefly review some part of height theory. 
We formulate it in terms of $\RR$-line bundles (not just line bundles), 
since $\RR$-line bundles are used in \cite{SiK3}. 
For details of height theory, we refer to \cite{La} and \cite{SiHt}.  

Let $K$ be a number field and $O_K$ its ring of integers.  
For $x = (x_0: \cdots: x_{n}) \in \PP^n(K)$, the logarithmic 
naive height of $x$ is defined by   
\[
h_{nv}(x) 
= 
\frac{1}{[K:\QQ]}
\left[ \sum_{P \in \Spec(O_K)\setminus\{0\}} 
\max_{0 \leq i \leq n} \{ -\ord_P(x_i) \} \log \#(O_K/P)
+ \sum_{\sigma: K \hookrightarrow \CC} 
\max_{0 \leq i \leq n} \{\log \vert\sigma(x_i)\vert \} \right]. 
\]
By the product formula, $h_{nv}(x)$ is independent of the homogeneous 
coordinates of $x$. It is also independent 
of the choice of fields over which $x$ is defined. 
Thus, we have the logarithmic naive height function 
$h_{nv}: \PP^n(\overline{\QQ}) \to \RR$. 

Let $X$ be a projective variety over $\overline{\QQ}$. 
Let $\Pic(X)$ denote the group of 
isomorphic classes of line bundles on $X$. 
For $L, M \in \Pic(X)$, 
we will often denote the tensor product 
additively, i.e., $L+M$ in place of $L\otimes M$.  
Let $F(X):= \{ f: X(\overline{\QQ}) \to \RR\}$ be   
the space of real-valued functions on $X(\overline{\QQ})$, 
and $B(X):= \{ f: X(\overline{\QQ}) \to \RR \mid 
\sup_{x \in X(\overline{\QQ})} |f(x)| < \infty\}$  
the space of real-valued bounded functions on $X(\overline{\QQ})$. 
We use the notation $h_1 = h_2 + O(1)$ if 
$h_1, h_2 \in F(X)$ satisfies $h_1 - h_2 \in B(X)$. 

\begin{Theorem}[cf. \cite{SiHt}]
\label{thm:height:machine}
For any projective variety $X$ over $\overline{\QQ}$, 
there exists a unique map 
\[
h_X : \Pic(X) \otimes_{\ZZ} \RR \longrightarrow 
\text{$F(X)$ modulo $B(X)$}, 
\quad L \mapsto h_{X,L}
\]
with the following properties\textup{:} 
\begin{enumerate}
\item[(1)] $h_{X, L}$ is $\RR$-linear;
\item[(2)] If $X = \PP^n$ and $L= \OO_{\PP^n}(1)$, then 
$h_{\PP^n, \OO_{\PP^n}(1)} = h_{nv} + O(1)$; 
\item[(3)] If $f: X \to Y$ is a morphism of projective varieties 
and $L$ is an element of $Pic(X) \otimes_{\ZZ} \RR$, then 
$h_{X, f^*L} = h_{Y, L}\circ f + O(1)$. 
\end{enumerate}
\end{Theorem}

\Proof
See \cite{SiHt},~\S3 for a proof when  
$\Pic(X) \otimes_{\ZZ} \RR$ is replaced by $\Pic(X)$, 
and $\RR$-linear by $\ZZ$-linear. 
(The theorem there holds without the assumption of the smoothness of $X$.) 
Let $h_X: \Pic(X) \to \text{$F(X)$ modulo $B(X)$}$ be the map defined 
in \cite{SiHt},~\S3. By tensoring $\RR$ over $\ZZ$, we have a map 
\[
h_X : \Pic(X) \otimes_{\ZZ} \RR \longrightarrow \text{$F(X)$ modulo $B(X)$}. 
\]
Obviously, $h_X$ satisfies (1) and (2). To see $h_X$ satisfies (3), 
let $L$ be an element of $\Pic(X) \otimes_{\ZZ} \RR$. We write 
$L \cong r_1 L_1 + \cdots + r_n L_n$ for $r_i \in \RR$ and 
$L_i \in \Pic(X)$.   
By \cite{SiHt},~\S3 it follows that $h_{X, f^*L_i} = h_{Y, L_i}\circ f$, 
and by (1) it follows that $h_{X, \sum_{i} r_i f^*L_i} 
= h_{Y, \sum_{i} r_i L_i}\circ f$. 
Thus we have $h_{X, f^*L} = h_{Y, L}\circ f$.  
The uniqueness of $h_X$ follows from \cite{SiHt},~\S3 
and $\RR$-linearity. 
\QED

By slight abuse of notation, we will use the same notation 
$h_{X,L} \in F(X)$ for a representative of 
$h_{X,L} \in \text{$F(X)$ modulo $B(X)$}$. We will often abbreviate  
$h_{X,L}$ by $h_L$. The following Northcott's finiteness theorem is 
useful in the study of rational points on projective varieties. 

\begin{Theorem}[Northcott's finiteness theorem]
\label{thm:Northcott}
Let $X$ be a projective variety over $\overline{\QQ}$, and
$L \in \Pic(X) \otimes_{\ZZ} \RR$ an ample $\RR$-line bundle.  
Then for any positive number $M$ and a positive integer $D$, 
the set
\[
\left\{
x \in X(\overline{\QQ}) \mid
[\QQ(x):\QQ] \leq D, \ h_L(x) \leq M 
\right\}
\]
is finite. 
\end{Theorem}

\Proof
See \cite{SiHt},~\S3 for a proof when 
$\Pic(X) \otimes_{\ZZ} \RR$ is replaced by $\Pic(X)$. 
In general, write $L \cong r_1 L_1 + \cdots + r_n L_n$ 
for $r_i > 0$ and ample $L_i$. 
By changing representatives if necessary, 
we can take $h_{L_i} \in F(X)$ such that 
$h_{L_i} \geq 0$ and $h_L = \sum_{i} r_i h_{L_i} \in F(X)$.
Then, the assertion follows from  
\[
\left\{
x \in X(\overline{\QQ}) \ \left\vert 
\begin{gathered}
\text{$[\QQ(x):\QQ] \leq D$,} \\
\text{$h_L(x) \leq M$}  
\end{gathered}
\right. \right\}
\quad \subseteq  \quad
\left\{
x \in X(\overline{\QQ}) \ \left\vert 
\begin{gathered}
\text{$[\QQ(x):\QQ] \leq D$,} \\
\text{$h_{L_1}(x) \leq \frac{M}{r_1}$}  
\end{gathered}
\right. \right\}.
\]
\QED

\subsection{Construction of canonical heights}
\label{subsec:construction:of:canonical:heights}
Let $X$ be a projective variety defined over a number field $K$, 
and $f_1, \cdots, f_k : X \to X$ be morphisms over $K$. 
Assume that there exists an $\RR$-line bundle $L$ on $X$  
such that $\sum_{i=1}^{k} {f_i^{*}(L)} 
\cong d L$ for some positive real 
number $d > k$ in $\Pic(X)\otimes\RR$. 

We set $\calF_0 := \{\id\}$ and  
$\calF_l := 
\{ f_{i_1}\circ\cdots\circ f_{i_l} \mid 1 \leq i_1,\cdots,i_l \leq k\}$ 
for $l \geq 1$. We set $\calF := \calF_1 (=\{f_1, \cdots, f_k\})$. 

\begin{Theorem}
\label{thm:main}
Let $X$ be a projective variety over a number field $K$, 
and $f_1, \cdots, f_k: X\to X$ morphisms over $K$. 
Assume that there exists an $\RR$-line bundle $L$ on $X$  
such that $\sum_{i=1}^{k} {f_i^{*}(L)} 
\cong {d} L$ for some positive real 
number $d > k$ in $\Pic(X)\otimes\RR$. 
Then there exists a unique real-valued function, 
\[
\widehat{h}_{L, \calF}: X(\overline{K}) \to \RR
\]
with the following properties\textup{:}
\begin{enumerate}
\item[(i)]
$\widehat{h}_{L, \calF} = h_L + O(1)$. 
\item[(ii)]
$\sum_{i=1}^{k} \widehat{h}_{L, \calF} (f_i(x)) = 
{d} \widehat{h}_{L, \calF}(x)$ 
for all $x \in X(\overline{K})$. 
We call $\widehat{h}_{L, \calF}$ the {\em canonical height function} 
associated with $L$ and $\calF$. 
\end{enumerate}
\end{Theorem}

\Proof
First let us construct $\widehat{h}_{L, \calF}$. 
We take any representative $h_L \in F(X)$ of 
$h_L \in \text{$F(X)$ modulo $B(X)$}$ and fix it. Since 
$\sum_{i=1}^{k} {f_i^{*}(L)} \cong {d} L$, 
there exists a constant $C > 0$ such that 
\[
\left\vert \sum_{i=1}^{k} h_L(f_i(x)) - {d} h_L(x) 
\right\vert \leq C
\]
for all $x \in X(\overline{K})$. 
Set $a_0(x):=h_L(x)$ and 
\[
a_l(x)
:= \frac{1}{{d}^l}\sum_{f \in \calF_l}h_L(f(x)) \in \RR.
\]
We claim $\{a_l(x)\}_{l=0}^{\infty}$ is a Cauchy sequence. 
Indeed, we have
\begin{align*}
\left\vert a_{l+1}(x) - a_l(x) \right\vert 
& =  
\left\vert
\frac{1}{{d}^{l+1}}\sum_{f \in \calF_{l+1}}h_L(f(x))
- \frac{1}{{d}^l}\sum_{f \in \calF_{l}}h_L(f(x))
\right\vert \\
& = 
\frac{1}{{d}^{l+1}}
\left\vert
\sum_{f \in \calF_{l}} \left\{
\sum_{i=1}^{k} h_L(f_i\circ f(x)) - {d} h_L(f(x))
\right\}\right\vert \\
& \leq 
\frac{\#\calF_{l}}{{d}^{l+1}} C 
= \frac{k^l C}{{d}^{l+1}}. 
\end{align*}
Since $\sum_{l=0}^{\infty} \frac{k^l C}{{d}^{l+1}} 
= \frac{C}{d-k} < \infty$, 
$\{a_l(x)\}_{l=0}^{\infty}$ is a Cauchy sequence. 

We define $\widehat{h}_{L, \calF}(x) := \lim_{l\to\infty} a_l(x)$.  
we need to show $\widehat{h}_{L, \calF}$ satisfies (i) and (ii). 
Note that $\left\vert a_{l}(x) - a_0(x) \right\vert 
\leq \sum_{\alpha=0}^{l-1} 
\left\vert a_{\alpha+1}(x) - a_{\alpha}(x)\right\vert
\leq \frac{C}{d-k}$. Since $a_0(x)=h_L(x)$,  
by letting $l\to\infty$ we obtain
\[
\left\vert \widehat{h}_{L, \calF}(x) -  h_L(x)\right\vert
\leq \frac{C}{d-k}.
\]
Since $\sum_{i=1}^{k} a_l(f_i(x)) = {d} a_{l+1}(x)$,
we also obtain (ii) by letting $l\to\infty$. 

To see the uniqueness of $\widehat{h}_{L, \calF}$, 
suppose both $\widehat{h}_1$ and $\widehat{h}_2$ satisfy 
(i) and (ii). By (i) there exists a constant $C'$ with 
$\vert \widehat{h}_1(x) - \widehat{h}_2(x)\vert \leq C'$ 
for all $x \in X(\overline{K})$. By (ii) it follows that 
$\widehat{h}_j(x) = \frac{1}{{d}^l}\sum_{f \in \calF_l} 
\widehat{h}_j(f(x))$. 
Since  
\[
\vert \widehat{h}_1(x) - \widehat{h}_2(x) \vert
\leq \frac{\#\calF_{l}}{{d}^{l}} C' 
= \frac{k^l C'}{{d}^{l}}, 
\]
we obtain $\widehat{h}_1(x) = \widehat{h}_2(x)$ 
by letting $l \to\infty$.
\QED

\subsection{Some properties of canonical heights}
\label{subsec:some:propertoes:of:canonical:heights}
Let the notation and assumption be as in \S1.2. 
For $x \in X(\overline{K})$, we define the forward
orbit of $x$ under $\calF$ to be 
\[
C(x):= \{f(x) \mid f \in \calF_l \,\,\text{for some $l \geq 0$}\}.   
\]
Note that when there is only one morphism $f_1$, i.e., $k=1$, 
$C(x)$ is finite if and only if $x$ is preperiodic with respect 
to $f_1$. 

\begin{Proposition}
\label{prop:property:of:height}
Let the notation and assumption be as in Theorem~\ref{thm:main}. 
Assume $L \in \Pic(X) \otimes_{\ZZ} \RR$ is ample. Then
\begin{enumerate}
\item[(1)]
$\widehat{h}_{L, \calF}(x) \geq 0$ for all $x \in X(\overline{K})$. 
\item[(2)]
$\widehat{h}_{L, \calF}(x) = 0$ if and only if 
$C(x)$ is finite. 
\end{enumerate}
\end{Proposition}

\Proof
(1): Take a height function $h_L \in F(X)$ corresponding to $L$. 
Since $L$ is ample, there exists a constant $C$ such that 
$h_L(x) \geq C$ for all $x \in X(\overline{K})$. Then we have 
\begin{align*}
\widehat{h}_{L, \calF}(x) & = 
\lim_{l \to\infty} \frac{1}{{d}^l}\sum_{f \in \calF_l}h_L(f(x)) \\
& \geq \limsup_{l \to\infty} \frac{1}{{d}^l} k^l C =0
\end{align*}  

(2): We first show the ``only if'' part.  
Let $K'$ be a finite extension of $K$ such that $x \in X(K')$. 
Since $\widehat{h}_{L, \calF} \geq 0$ and 
$\sum_{i=1}^{k} \widehat{h}_{L, \calF} (f_i(x)) = 
{d} \widehat{h}_{L, \calF}(x)$, 
it follows that $\widehat{h}_{L, \calF} (f(x)) = 0$ for 
all $f \in \bigcup_{l\geq 0}\calF_l$. Then by Northcott's finiteness theorem 
(Theorem~\ref{thm:Northcott}),  
$C(x) \subseteq \{y \in X(K') \mid \widehat{h}_{L, \calF} (y) =0\}$ 
is finite. 

We next show the ``if'' part. 
Suppose $a:=\widehat{h}_{L, \calF}(x) > 0$. 
Since $\sum_{i=1}^{k} \widehat{h}_{L, \calF} (f_i(x)) = 
{d} \widehat{h}_{L, \calF}(x)$, it follows that 
$\widehat{h}_{L, \calF} (f_{i_1}(x)) \geq 
\frac{d}{k} a$ for some $i_1$. 
Since $\sum_{i=1}^{k} \widehat{h}_{L, \calF} (f_i(f_{i_1}(x))) = 
{d} \widehat{h}_{L, \calF}(f_{i_1}(x))$, it follows that 
$\widehat{h}_{L, \calF} (f_{i_2}\circ f_{i_1}(x)) \geq 
(\frac{d}{k})^2 a$ for some $i_2$. 
Similarly, for any $l \in \ZZ_{> 0}$, there exists $f \in \calF_l$ 
with $\widehat{h}_{L, \calF}(f(x)) \geq (\frac{d}{k})^l a$. 
Since $\frac{d}{k} > 1$, $C(x)$ cannot be finite. 
\QED

\begin{Corollary}
\label{cor:property:of:height}
Let the notation and assumption be as in Theorem~\ref{thm:main}. 
Assume $L \in \Pic(X) \otimes_{\ZZ} \RR$ is ample. Then
for any $D \in \ZZ_{> 0}$, 
\begin{equation*}
\left\{
x \in X(\overline{K}) \ \left\vert \ 
\text{$[K(x): K] \leq D$, \ $C(x)$ is finite} 
\right. \right\}
\end{equation*}
is finite. 
\end{Corollary}

\Proof
The assertion follows from Proposition~\ref{prop:property:of:height}(2) 
and Theorem~\ref{thm:Northcott}.
\QED

\subsection{Examples}
\label{subsec:examples}
\subsubsection{Abelian varieties}
\label{subsec:abelian:var}
Let $A$ be an abelian variety over $\overline{\QQ}$, and 
$L$ a symmetric ample line bundle on $A$.  
Since $[2]^*(L)\cong 4L$, the canonical height function 
$\widehat{h}_{L, [2]}:A(\overline{\QQ}) \to \RR_{\geq 0}$, called 
N\'eron-Tate's height function, is defined so 
that $\widehat{h}_{L, [2]}$ is a height corresponding to $L$ 
and $\widehat{h}_L\circ [2] = 4 \widehat{h}_L$. 

\subsubsection{Projective spaces}
\label{subsec:proj:space}
Let $F_0, \cdots, F_N \in \overline{\QQ}[X_0, \cdots, X_N]$ be 
homogeneous polynomials of degree $d \geq 1$ such that $0$ is the only  
common zero of $F_0, \cdots, F_N$. Then $f=(F_0: \cdots: F_N)$ defines 
a morphism of $\PP^N$ to $\PP^N$. 
Assume $d \geq 2$. 
Since $f^*\OO_{\PP^N}(1) \cong \OO_{\PP^N}(d)$, 
the canonical height function $\widehat{h}_{\OO_{\PP^N}(1), f}$ 
associated with $\OO_{\PP^N}(1)$ and $f$ is defined. 

More generally, for $i=1, \cdots, k$, let 
$f_i: \PP^N \to \PP^N$ be morphisms
of degree $d_i \geq 1$. Assume $d_i \geq 2$ for some $i$. 
Since $\bigotimes_{i=1}^{k} f_i^{*}(\OO_{\PP^N}(1)) \simeq
\OO_{\PP^N}(\sum_{i=1}^k d_i)$, 
the canonical height function $\widehat{h}_{\OO_{\PP^N}(1), 
\{f_1, \cdots, f_k\}}$ associated with 
$\OO_{\PP^N}(1)$ and $\{f_1, \cdots, f_k\}$ is defined. 
We note that the dynamical system $(\PP^N; f_1, \cdots, f_k)$ 
over $\CC$ in the case $N=1$ (dynamics of finitely generated 
rational semigroups) is studied in \cite{HM}. 

\subsubsection{Toric varieties}
\label{subsec:toric:varieties}
Let $\PP(\Delta)$ be a smooth projective toric variety 
$\overline{\QQ}$, and $D$ an ample 
torus-invariant Cartier divisor on $\PP(\Delta)$.  
Let $p \geq 2$ be an integer, and 
$[p]: \PP(\Delta) \to \PP(\Delta)$ the morphism associated 
with the multiplication of $\Delta$ by $p$. 
Since $[p]^*(\OO_{\PP(\Delta)}(D)) \simeq \OO_{\PP(\Delta)}(D)^{\otimes p}$, 
the canonical height function 
$\widehat{h}_{\OO_{\PP(\Delta)}(D), [p]}: \PP(\Delta)(\overline{\QQ}) \to \RR$ 
associated with $\OO_{\PP(\Delta)}(D)$ and $[p]$ is defined. 
For details, we refer Maillot \cite{Mai},~Theorem~3.3.3. 

\subsubsection{K3 surfaces, I}
\label{subsubsec:K3:I}
Let $S$ be a K3 surface in $\PP^2\times\PP^2$ given by the
intersection of two hypersurfaces of bidegrees $(1,1)$ and $(2,2)$ 
over $\overline{\QQ}$. 
Such K3 surfaces are called Wheler's K3 surfaces (cf. \cite{Maz}).  
Silverman \cite{SiK3} constructed a canonical height function 
on $S$ and developed an arithmetic theory. Here, we will show 
we can construct a height function on $S$ by using 
Theorem~\ref{thm:main}, which turns out to coincide with 
the one Silverman constructed up to a constant multiple. 

Let $p_i: S\to\PP^2$ be the projection to the $i$-th factor 
for $i=1,2$. Since $p_i$ is a double cover, it gives an 
involution $\sigma_i \in \Aut(S)$. Set $L_i := p_i^*\OO_{\PP^2}(1)$. 
By \cite{SiK3},~Lemma~2.1, we have 
\begin{align*}
\sigma_1^*(L_1) \cong L_1, \quad 
\sigma_2^*(L_1) \cong - L_1 +4L_2,\\
\sigma_1^*(L_2) \cong 4 L_1 -L_2, \quad 
\sigma_2^*(L_2) \cong L_2.
\end{align*}
Thus the ample line bundle $L:=L_1+L_2$ on $S$ satisfies 
$\sigma_1^*(L) + \sigma_2^*(L) \cong 4 L$. 
Applying Theorem~\ref{thm:main} to $S$, $L$ and $\{\sigma_1, \sigma_2\}$, 
we obtain a height function 
$\widehat{h}_{L, \{\sigma_1,\sigma_2\}}: S(\overline{\QQ}) \to \RR_{\geq 0}$.

Let us recall Silverman's construction of a height function. 
Setting $E^{+}:= (2+\sqrt{3})L_1 - L_2$ and $E^{-}:= -L_1 +(2+\sqrt{3})L_2$,  
we have $(\sigma_2\circ\sigma_1)^{\pm 1 *}(E^{\pm}) 
= (7+4\sqrt{3}) E^{\pm}$. 
Then we obtain two height functions $\widehat{h}^{\pm}$ 
which are characterized by being 
a height corresponding to $E^{\pm}$ and 
$\widehat{h}^{\pm} \circ (\sigma_2\circ\sigma_1)^{\pm 1}  = (7+4\sqrt{3}) 
\widehat{h}^{\pm}$. (In the notation of Theorem~\ref{thm:main}, 
$\widehat{h}^{+} = \widehat{h}_{E^{+}, {\sigma_2\circ\sigma_1}}$ and 
$\widehat{h}^{-} = \widehat{h}_{E^{-}, {\sigma_1\circ\sigma_2}}$.)
Silverman defined a canonical height function 
on $S$ to be $\widehat{h}^{+} + \widehat{h}^{-}$, which is
denoted by $\widehat{h}_{Sil}$ in the following. 

The next proposition relates 
$\widehat{h}_{Sil}$ with $\widehat{h}_{L, \{\sigma_1,\sigma_2\}}$. 

\begin{Proposition}
\label{prop:h:and:Sil}
$\widehat{h}_{Sil} = (1+\sqrt{3}) \ 
\widehat{h}_{L, \{\sigma_1,\sigma_2\}}$.
\end{Proposition}

\Proof
It suffices to show $\frac{1}{1+\sqrt{3}}\widehat{h}_{Sil}$ 
satisfies the two conditions of Theorem~\ref{thm:main}. 
Since $\widehat{h}_{Sil}$ is by construction a height 
corresponding to $E^{+} + E^{-} \cong (1+\sqrt{3}) L$, 
it follows that $\frac{1}{1+\sqrt{3}}\widehat{h}_{Sil}$ satisfies 
the first condition. 
By \cite{SiK3},~Theorem~1.1(ii), we have 
\[
\widehat{h}^{\pm}\circ\sigma_1 = (2+\sqrt{3})^{\mp 1} 
\widehat{h}^{\mp}, \quad
\widehat{h}^{\pm}\circ\sigma_2 = (2+\sqrt{3})^{\pm 1} 
\widehat{h}^{\mp}. 
\]
Then we have 
\begin{align*}
& \frac{1}{1+\sqrt{3}}\widehat{h}_{Sil}(\sigma_1(x)) + 
  \frac{1}{1+\sqrt{3}}\widehat{h}_{Sil}(\sigma_2(x)) \\
= & \frac{1}{1+\sqrt{3}} \left\{
\widehat{h}^{+}(\sigma_1(x))+ \widehat{h}^{-}(\sigma_1(x))
+\widehat{h}^{+}(\sigma_2(x))+ \widehat{h}^{-}(\sigma_2(x)) \right\}\\
= & \frac{1}{1+\sqrt{3}}
\left((2+\sqrt{3})+ (2+\sqrt{3})^{-1}\right)
(\widehat{h}^{+}(x) + \widehat{h}^{-}(x))
= 4 \frac{1}{1+\sqrt{3}}\widehat{h}_{Sil}.
\end{align*}
Thus $\frac{1}{1+\sqrt{3}}\widehat{h}_{Sil}$ also satisfies the 
second condition of Theorem~\ref{thm:main}. 
\QED

\subsubsection{K3 surfaces, II}
Let $S$ be a K3 surface in $\PP^2\times\PP^2$ given by the
intersection of two hypersurfaces of bidegrees $(1,2)$ and $(2,1)$ 
over $\overline{\QQ}$. 
Let $p_i: S\to\PP^2$ be the projection to the $i$-th factor 
for $i=1,2$. Since $p_i$ is a double cover, it gives an 
involution $\sigma_i \in \Aut(S)$. Set $L_i := p_i^*\OO_{\PP^2}(1)$. 
By similar computations as above, we have 
\begin{gather*}
\sigma_1^*(L_1) \cong L_1, \quad 
\sigma_2^*(L_1) \cong - L_1 +5L_2, \\
\sigma_1^*(L_2) \cong 5 L_1 -L_2, \quad 
\sigma_2^*(L_2) \cong L_2.
\end{gather*}
Thus the ample line bundle $L:=L_1+L_2$ on $S$ satisfies 
$\sigma_1^*(L) + \sigma_2^*(L) \cong 5 L$. 
Applying Theorem~\ref{thm:main} to $S$, $L$ and $\{\sigma_1, \sigma_2\}$, 
we obtain a height function 
$\widehat{h}_{L, \{\sigma_1,\sigma_2\}}: S(\overline{\QQ}) \to \RR_{\geq 0}$. 

\subsubsection{K3 surfaces, III}
\label{subsec:K3:surfaces:III}
Let $S$ be a hypersurface of tridegrees $(2,2,2)$ in 
$\PP^1\times\PP^1\times\PP^1$ over $\overline{\QQ}$ (cf. \cite{Maz}). 
For $i=1,2,3$, let $p_i: S\to\PP^1\times\PP^1$ be the projection to 
the $(j,k)$-th factor with $\{i,j,k\}=\{1,2,3\}$. 
Since $p_i$ is a double cover, it gives an 
involution $\sigma_i \in \Aut(S)$. 
Let $q_i: S \to \PP^1$ be the projection to the $i$-th factor, 
and  set  $L_i := q_i^*\OO_{\PP^1}(1)$. 
By similar computations as above, we have 
\begin{gather*}
\sigma_i^*(L_i) \cong -L_i + 2L_j +2L_k 
\quad\text{for $\{i,j,k\}=\{1,2,3\}$}, \\ 
\sigma_j^*(L_i) \cong L_i \quad\text{for $i\neq j$}.
\end{gather*}
Thus the ample line bundle $L:=L_1+L_2+L_3$ on $S$ satisfies 
$\sigma_1^*(L) + \sigma_2^*(L) + \sigma_3^*(L) \cong 5 L$. 
Applying Theorem~\ref{thm:main} to $S$, $L$ and $\{\sigma_1, \sigma_2, 
\sigma_3\}$, we obtain a height function $
\widehat{h}_{L, \{\sigma_1,\sigma_2,\sigma_3\}}: S(\overline{\QQ}) 
\to \RR_{\geq 0}$.
We remark that Wang \cite{Wa} constructed a height function on $S$. 
It may be interesting to compare it with 
$\widehat{h}_{L, \{\sigma_1, \sigma_2, \sigma_3\}}$.  

To illustrate Proposition~\ref{prop:property:of:height}(2), 
suppose now $S$ is defined by the affine equation 
\[
x(1-x)+y(1-y)+z(1-z) - xyz = 0.  
\]
Then $S$ is smooth. It is easy to check that the orbit of $(0,0,0)$ 
under $\{\sigma_1,\sigma_2, \sigma_3\}$ is 
\[
\{(0,0,0),(1,0,0),(0,1,0),(0,0,1),(1,1,0),(1,0,1),(0,1,1)\}
\] 
and thus finite. In particular, 
$\widehat{h}_{L, \{\sigma_1,\sigma_2,\sigma_3\}}((0,0,0))=0$. 

\subsubsection{Product of a projective space and an abelian variety}
Let $A$ be an abelian variety over $\overline{\QQ}$, and 
$L$ a symmetric ample line bundle on $A$.  
Let $f=(F_0: \cdots: F_N): \PP^N \to \PP^N$ be a morphism 
defined by the homogeneous polynomials $F_0, \cdots, F_N$ 
of degree $d > 1$ such that $0$ is the only common zero of 
$F_0, \cdots, F_N$. 
Set $X= A \times \PP^N$, $g_1=[2]\times \id_{\PP^N}$, and 
$g_2=\id_{A} \times f$. 
Put $M = p_1^* L \otimes p_2^* \OO_{\PP^N}(1)$, 
where $p_1$ and $p_2$ are the obvious projections. 
Then 
\[
\overbrace{g_{1}^*(M) \otimes \cdots \otimes g_{1}^*(M)
}^{\text{$(d-1)$ times}} \otimes 
g_2^*(M) \otimes g_2^*(M) \otimes g_2^*(M) 
\simeq M^{\otimes (4d-1)}.
\]
Thus we have the canonical height function 
$\widehat{h}_{M, \{g_1,\cdots, g_1, g_2, g_2, g_2\}}$ on 
$A \times \PP^N(\overline{\QQ})$. 
It is easy to see 
$\widehat{h}_{M, \{g_1,\cdots, g_1, g_2, g_2, g_2\}} 
= \widehat{h}_{L} \circ p_1 
+ \widehat{h}_{\OO_{\PP^N}(1), f} \circ p_2$. 

\medskip
\section{Canonical heights of subvarieties}
\label{sec:height:subvar}
In this section, first we briefly review the adelic intersection 
theory due to Zhang \cite{Zh}. Then we construct 
an adelic sequence for a given dynamical system of $k$ morphisms 
associated with a line bundle, 
and finally define canonical heights for its subvarieties.  
We refer to \cite{Zh} and \cite{Mo} for details of adelic 
intersection theory, to \cite{So} and \cite{BGS} 
for details of Gillet-Soul\'e's arithmetic intersection 
theory and height theory of subvarieties in general. 

\subsection{Adelic sequence}
Let $K$ be a number field and $O_K$ its ring of integers. 
By a projective arithmetic variety, we mean 
a flat and projective integral scheme over $O_K$. 
A pair $\overline{\calL}=(\calL, \Vert\cdot\Vert)$ is called 
a {\em $C^{\infty}$-hermitian $\QQ$-line bundle} over 
a projective arithmetic variety $\calX$ if  
$\calL$ is a $\QQ$-line bundle over $\calX$ 
and $\Vert\cdot\Vert$ is a $C^{\infty}$-hermitian metric 
on $\calL \otimes_{\ZZ} \CC$ that is invariant under the complex 
conjugation. Recall that $\Vert\cdot\Vert$ is said to be $C^{\infty}$ 
if, for any analytic morphism $h: M \to \calX \otimes_{\ZZ}{\CC}$ 
from any complex manifold $M$, $h^{*}(\Vert\cdot\Vert)$ 
is a $C^{\infty}$-hermitian metric on $h^*(\calL \otimes_{\ZZ}{\CC})$. 

Following \cite{Mo1}, we say that 
$\overline{\calL}$ is {\em nef} 
if the first Chern form $\chern_1(\calL \otimes_{\ZZ} \CC, 
\Vert\cdot\Vert)$ is semipositive on $X(\CC)$, 
and for all one-dimensional integral closed subschemes $\Gamma$ 
of $\calX$, $\adeg(\rest{\overline{\calL}}{\Gamma}) \geq 0$. 
We say that $\overline{\calL}$ is {\em $\QQ$-effective} 
(denoted by $\overline{\calL} \gtrsim 0$) if 
there is a positive integer $n$ and a non-zero section 
$s \in H^0(\calX, \calL^{\otimes n})$ with 
$\Vert s \Vert_{\sup} \leq 1$. 
For two $C^{\infty}$-hermitian $\QQ$-line bundles 
$\overline{\calL}$ and ${\overline{\calM}}$, 
we use the notation $\overline{\calL} \gtrsim \overline{\calM}$ 
if $\overline{\calL} \otimes \overline{\calM}^{-1}
\gtrsim 0$. We quote \cite{Mo1},~Proposition~2.3 here. 

\begin{Proposition}[\cite{Mo1}]
\label{prop:nef:Q:effective}
Set $\dim X = e$. 
\renewcommand{\labelenumi}{(\arabic{enumi})}
\begin{enumerate}
\item
If $\overline{\calL_1}, \cdots,\overline{\calL_{e + 1}}$ 
are nef, 
then 
$\adeg(\achern_1(\overline{\calL_1})\cdots 
\achern_1(\overline{\calL_{e + 1}})) \geq 0$. 
\item
If $\overline{\calL_1}, \cdots,\overline{\calL_{e}}$ 
are nef and $\overline{\calM}$ is $\QQ$-effective, 
then 
$\adeg(\achern_1(\overline{\calL_1})\cdots 
\achern_1(\overline{\calL_{e}})\cdot
\achern_1(\overline{\calM})) \geq 0$. 
\end{enumerate}
\end{Proposition}

Let $X$ a projective variety over $K$ and 
$L$ a $\QQ$-line bundle over $X$. 
A pair $(\calX, \overline{\calL})$ of a projective 
arithmetic variety and a $C^{\infty}$-hermitian $\QQ$-line bundle 
is called a {\em $C^{\infty}$-model of $(X,L)$} 
if one has $\calX\otimes_{O_K}K=X$ and $\calL\otimes_{O_K}K=L$. 

Following \cite{Zh} and \cite{Mo}, we make following definitions. 
\begin{Definition}
\label{def:adelic:sequence}
\renewcommand{\labelenumi}{(\arabic{enumi})}
\begin{enumerate}
\item
Let $\{(\calX_n, \overline{\calL_n})\}_{n=0}^{\infty}$ be 
a sequence of $C^{\infty}$-models of $(X,L)$. We say that 
$\{(\calX_n, \overline{\calL_n})\}_{n=0}^{\infty}$ is an
{\em adelic sequence} of $(X,L)$ if it satisfies the following 
three conditions:
\renewcommand{\labelenumii}{(\roman{enumii})}
\begin{enumerate}
\item
For every $n \geq 0$, $\overline{\calL_n}$ is nef;
\item
There exists a Zariski dense open set $U \subseteq \Spec(O_K)$ 
such that, for every $n \geq 0$, $\rest{\calX_n}{U} 
= \rest{\calX_{n+1}}{U}$ and $\rest{\calL_n}{U} = \rest{\calL_{n+1}}{U}$; 
We set $\Spec{O_K}\setminus U = \{\frakP_1, \cdots, \frakP_r\}$; 
\item
For every $n \geq 0$, there exist 
a projective arithmetic variety $\calX_n^{'}$, 
birational morphisms $\mu_n: \calX_n^{'} \to \calX_n$ and 
$\nu_n: \calX_n^{'} \to \calX_{n+1}$, and positive rational numbers 
$c_{n1}, \cdots, c_{nr} \in \QQ_{>0}$ and 
positive real numbers $c_{n\sigma} \in \RR_{>0}$ for every 
$\sigma: K \hookrightarrow \CC$ 
such that 
\begin{multline*}
\pi_{\calX_n^{'}}^*\left(\left(
-\sum_{\alpha=1}^r c_{n\alpha} [\frakP_{\alpha}], 
-\sum_{\sigma:K\hookrightarrow \CC} c_{n\sigma} [\sigma]
\right)\right) \lesssim
\mu_n^*(\overline{\calL_n}) - \nu_n^*(\overline{\calL_{n+1}}) \\
\lesssim
\pi_{\calX_n^{'}}^*\left(\left(  \sum_{\alpha=1}^r c_{n\alpha} 
[\frakP_{\alpha}]
\sum_{\sigma:K\hookrightarrow \CC} c_{n\sigma} [\sigma]
\right)\right),
\end{multline*}
\[
\sum_{n=0}^{\infty} c_{n\alpha} < \infty \quad (\alpha = 1, \cdots, r), 
\qquad
\sum_{n=0}^{\infty} c_{n\sigma} < \infty \quad (\sigma:K\hookrightarrow\CC)
\]
where $\pi_{\calX_n^{'}}: \calX_n^{'} \to \Spec(O_K)$ is the structure 
morphism. 
\end{enumerate}
A open set $U$ of $\Spec(O_K)$ satisfying (ii) and (iii) is 
called a common base of $\{(\calX_n, \overline{\calL_n})\}_{n=0}^{\infty}$. 
Note that if $U$ is a common base of 
$\{(\calX_n, \overline{\calL_n})\}_{n=0}^{\infty}$ and $U'$ is 
a non-empty Zariski open set of $U$, then $U'$ is also a common base 
of $\{(\calX_n, \overline{\calL_n})\}_{n=0}^{\infty}$. 
\item
Two adelic sequences 
$\{(\calX_n, \overline{\calL_n})\}_{n=0}^{\infty}$ 
and $\{(\calY_n, \overline{\calM_n})\}_{n=0}^{\infty}$
of $(X,L)$ are said to be {\em equivalent} if 
there exists a non-empty Zariski open set $\widetilde{U}$ of 
$\Spec(O_K)$ with the following properties: 
\renewcommand{\labelenumii}{(\roman{enumii})}
\begin{enumerate}
\item
$\widetilde{U}$ is a common base of both 
$\{(\calX_n, \overline{\calL_n})\}_{n=0}^{\infty}$ and 
$\{(\calY_n, \overline{\calM_n})\}_{n=0}^{\infty}$. 
\item 
For every $n \geq 0$, there exist 
a projective arithmetic variety $\calZ_n$, 
birational morphisms $\widetilde{\mu}_n: \calZ_n \to \calX_n$ and 
$\widetilde{\nu}_n: \calZ_n \to \calY_n$, 
and positive rational numbers 
$\widetilde{c}_{n1}, \cdots, \widetilde{c}_{n\widetilde{r}}$ 
and positive real numbers $\widetilde{c}_{n\sigma}$ 
for every $\sigma:K\hookrightarrow\CC$ such that 
\begin{multline*}
\pi_{\calZ_n}^{*} 
\left(\left( - \sum_{\alpha=1}^{\widetilde{r}} {\widetilde{c}}_{n\alpha} 
[{\widetilde{\frakP}}_{\alpha}], 
-\sum_{\sigma:K\hookrightarrow \CC} \widetilde{c}_{n\sigma} [\sigma]
\right)\right) 
\lesssim
\widetilde{\mu}_n^*(\overline{\calL_n}) 
- \widetilde{\nu}_n^*(\overline{\calM_{n}}) \\ 
\lesssim
\pi_{\calZ_n^{'}}^*
\left(\left(  \sum_{\alpha=1}^{\widetilde{r}} {\widetilde{c}}_{n\alpha} 
[{\widetilde{\frakP}}_{\alpha}],
\sum_{\sigma:K\hookrightarrow \CC} \widetilde{c}_{n\sigma} [\sigma]
\right)\right), 
\end{multline*}
\[
\sum_{n=0}^{\infty} {\widetilde{c}}_{n\alpha} < \infty \quad (\alpha = 1, \cdots, {\widetilde{r}}), 
\qquad 
\sum_{n=0}^{\infty} \widetilde{c}_{n\sigma} < \infty \quad (\sigma:K\hookrightarrow\CC)
\]
where $\Spec{O_K}\setminus \widetilde{U} = 
\{\widetilde{\frakP}_1, \cdots, \widetilde{\frakP}_{\widetilde{r}}\}$ 
and $\pi_{\calZ_n}: \calZ_n \to \Spec(O_K)$ is the structure 
morphism. 
\end{enumerate}
\end{enumerate}
\end{Definition}

The next lemma follows directly from Definition~\ref{def:adelic:sequence}. 

\begin{Lemma}
\label{lemma:adelic:sequence:associated}
Let $\{(\calX_n, \overline{\calL_n})\}_{n=0}^{\infty}$ be an 
adelic sequence of $(X,L)$. 
\renewcommand{\labelenumi}{(\arabic{enumi})}
\begin{enumerate}
\item
Let $Y$ be a subvariety of $X_{\overline{K}}$ 
such that $Y$ is defined over $K$. 
Let $\calY_n$ be the Zariski closure of $Y$ in $\calX_n$. 
Then $\{(\calY_n, \rest{\overline{\calL_n}}{\calY_n})\}_{n=0}^{\infty}$ 
is an adelic sequence of $(Y, \rest{L}{Y})$. 
\item Let $K'$ be a finite extension field of $K$ 
and $O_{K'}$ its ring of integers. 
We set $\calX^{'}_n = \calX_n \times_{\Spec(O_K)} \Spec(O_K')$ 
and $\overline{\calL^{'}_n} = \overline{\calL^{'}_n} 
\otimes_{O_K} O_K'$. Then 
$\{(\calX^{'}_n, \overline{\calL{'}_n})\}_{n=0}^{\infty}$ is an 
adelic sequence of $(X_{K'} ,L_{K'})$. 
\end{enumerate}
\end{Lemma}

The next theorem (a special case of \cite{Zh},~Theorem~(1.4)
and \cite{Mo},~Proposition~4.1.1) asserts the existence 
of the adelic intersection number. 

\begin{Theorem}[A special case of adelic intersection number \cite{Zh}]
\label{thm:adelic:intersection}
Let $\{(\calX_n, \overline{\calL_n})\}_{n=0}^{\infty}$ be an 
adelic sequence of $(X,L)$. Then the arithmetic 
intersection number 
\[
\adeg\left(\achern_1(\overline{\calL_n})^{\dim X +1} \right)
\]
converges as $n$ goes to $\infty$. 
Moreover, the limit is independent of the choice of equivalent 
adelic sequences. We call the limit the adelic intersection number 
of $\{(\calX_n, \overline{\calL_n})\}_{n=0}^{\infty}$. 
\end{Theorem}

\Proof
Since we adopt a slightly different definition for 
an adelic sequence, we sketch a proof here.   
Let $U$ be a common base and 
$\mu_n: \calX_n^{'} \to \calX_n$ and 
$\nu_n: \calX_n^{'} \to \calX_{n+1}$ 
birational morphisms as in Definition~\ref{def:adelic:sequence}. 
Then 
\begin{multline*}
\achern_1(\overline{\calL_n})^{\dim X +1}
- \achern_1(\overline{\calL_{n+1}})^{\dim X+1} \\
= 
\achern_1(\mu_n^*(\overline{\calL_n}))^{\dim X +1}
- \achern_1(\nu_n^*(\overline{\calL_{n+1}}))^{\dim X+1} \\
= 
\achern_1\left(\mu_n^*(\overline{\calL_n}) 
- \nu_n^*(\overline{\calL_{n+1}})\right)
\left\{
\sum_{i=0}^{\dim X} 
\achern_1(\mu_n^*(\overline{\calL_n}))^i 
\achern_1(\nu_n^*(\overline{\calL_n}))^{\dim X -i} 
\right\},
\end{multline*}
where we use the projection formula in the first equality. 
Since 
\begin{multline*}
\pi_{\calX_n^{'}}^*\left(\left( - \sum_{\alpha=1}^r c_{n\alpha} 
[\frakP_{\alpha}], 
- \sum_{\sigma:K\hookrightarrow \CC} c_{n\sigma} [\sigma]
\right)\right) \lesssim
\mu_n^*(\overline{\calL_n}) 
- \nu_n^*(\overline{\calL_{n+1}}) \\
\lesssim
\pi_{\calX_n^{'}}^*\left(\left( \sum_{\alpha=1}^r c_{n\alpha} 
[\frakP_{\alpha}], 
\sum_{\sigma:K\hookrightarrow \CC} c_{n\sigma} [\sigma]
\right)\right), 
\end{multline*}
it follows from the projection formula and 
Proposition~\ref{prop:nef:Q:effective} that 
\begin{align*}
& \left\vert
\achern_1(\mu_n^*(\overline{\calL_n}))^{\dim X +1}
- \achern_1(\nu_n^*(\overline{\calL_{n+1}}))^{\dim X+1}
\right\vert \\
& \leq
\pi_{\calX_n}^*\left(\left(\sum_{\alpha=1}^r c_{n\alpha} 
[\frakP_{\alpha}], 
\sum_{\sigma:K\hookrightarrow \CC} c_{n\sigma} [\sigma] 
\right)\right) 
\left\{
\sum_{i=0}^{\dim X} 
\achern_1(\mu_n^*(\overline{\calL_n}))^i 
\achern_1(\nu_n^*(\overline{\calL_n}))^{\dim X -i} 
\right\}
\\
& = \left(\sum_{\alpha=1}^r  c_{n\alpha} \log\#(O_K/\frakP_{\alpha})
+ \frac{1}{2} \sum_{\sigma:K\hookrightarrow \CC} c_{n\sigma}
\right) \chern_1(L)^{\dim X}.   
\end{align*}
Since $\sum_{n=0}^{\infty} c_{n\alpha} < \infty$ 
and $\sum_{n=0}^{\infty} c_{n\sigma} < \infty$, 
$\adeg\left(\achern_1(\overline{\calL_n})^{\dim X +1} \right)$ 
converges as $n$ goes to $\infty$. 

Next let $\{(\calX_n, \overline{\calL_n})\}_{n=0}^{\infty}$ and 
$\{(\calY_n, \overline{\calM_n})\}_{n=0}^{\infty}$ be equivalent 
adelic sequences for $(X, L)$. Take birational morphisms 
$\widetilde{\mu}_n: \calZ_n \to \calX_n$ and 
$\widetilde{\nu}_n: \calZ_n \to \calY_n$ 
as in Definition~\ref{def:adelic:sequence}. 
Note that by the projection formula 
$\achern_1(\overline{\calL_n})^{\dim X +1} 
= \achern_1(\widetilde{\mu}_n^*(\overline{\calL_n}))^{\dim X +1}$ 
and 
$\achern_1(\overline{\calM_n})^{\dim X +1} 
= \achern_1(\widetilde{\nu}_n^*(\overline{\calM_n}))^{\dim X +1}$. 
A similar argument as above yields 
\[
\left\vert
\achern_1(\widetilde{\mu}_n^*(\overline{\calL_n}))^{\dim X +1}
- \achern_1(\widetilde{\nu}_n^*(\overline{\calM_{n}}))^{\dim X+1}
\right\vert
\leq
\left(\sum_{\alpha=1}^{\widetilde{r}}  
\widetilde{c}_{n\alpha} \log\#(O_K/\widetilde{\frakP}_{\alpha})
+ \frac{1}{2} \sum_{\sigma:K\hookrightarrow \CC} \widetilde{c}_{n\sigma}
\right) \chern_1(L)^{\dim X}. 
\]
Thus $\adeg\left(\achern_1(\overline{\calL_n})^{\dim X +1} \right)$ and 
$\adeg\left(\achern_1(\overline{\calM_n})^{\dim X +1} \right)$ 
have the same limit as $n$ goes to $\infty$.
\QED

\subsection{Adelic sequence arising from several morphisms}
\label{subsec:adelic:sequence:from:morphisms}
Let $(X; f_1, \cdots, f_k)$ be a dynamical system 
of $k$ morphisms over $K$ associated with $L$ of degree $d > k$.  
{\em In \S\ref{subsec:adelic:sequence:from:morphisms} we assume that 
$X$ is normal, $f_1, \cdots f_k$ are surjective 
and $L$ is an ample $\QQ$-line bundle. }

Since $X$ is normal and $L$ is ample, we can take a 
$C^{\infty}$-model $(\calX, \overline{\calL})$ of $(X,L)$ 
such that $\calX$ is normal and $\overline{\calL}$ is nef. 
We take a non-empty Zariski open subset $U \subseteq \Spec(O_K)$ 
such that each $f_i: X \to X$ extends to $f_i: \calX_{U} \to \calX_{U}$ 
and that $f_1^*(\calL_{U})\otimes \cdots 
\otimes f_k^*(\calL_{U}) \simeq \calL_{U}^{\otimes d}$ 
in $\Pic(\calX_U)\otimes\QQ$. By fixing an isomorphism, 
we use the notation $ 
f_1^*(\calL_{U})\otimes \cdots 
\otimes f_k^*(\calL_{U}) = \calL_{U}^{\otimes d}$ in the following.  

Let us construct another $C^{\infty}$-model of $(X,L)$ from 
$(\calX, \overline{\calL})$. 
Let $\calX^{(i)}$ be the normalization of 
\[
\calX_U \overset{f_i}{\longrightarrow} \calX_U \hookrightarrow \calX.
\]
We denote the induced morphism by $f^{(i)}: \calX^{(i)} \to \calX$. 
Let $\calX_1$ be the Zariski closure of 
\[
\calX_U \overset{\Delta}{\longrightarrow}
\overbrace{\calX_U\times_U\cdots\times_U\calX_U}^{k}
\hookrightarrow 
\calX^{(1)}\times_{O_K}\cdots\times_{O_K}\calX^{(k)}, 
\]
where $\Delta$ is the diagonal map. 
Let $p_i: \calX^{(1)}\times_{O_K}\cdots\times_{O_K}\calX^{(k)} 
\to \calX^{(i)}$ denote the projection to the $i$-th factor. 
We set 
\[
\overline{\calL_1} = \left[\rest{\left((f^{(1)}\circ p_1)^*\overline{\calL}
\otimes\cdots\otimes (f^{(k)}\circ p_k)^*\overline{\calL}\right)}{\calX_1}
\right]^{\otimes \frac{1}{d}}.
\]
Then we get the $C^{\infty}$-model $(\calX_1, \overline{\calL_1})$ 
of $(X, L)$ from $(\calX_0, \overline{\calL_0}) := 
(\calX, \overline{\calL})$. 
Inductively we obtain the $C^{\infty}$-model 
$(\calX_{n+1}, \overline{\calL_{n+1}})$ of $(X,L)$ 
from $(\calX_n, \overline{\calL_n})$. 

\begin{Theorem}
\label{thm:adelic:sequence:from:morphisms}
\renewcommand{\labelenumi}{(\arabic{enumi})}
\begin{enumerate}
\item
$\{(\calX_n, \overline{\calL_n})\}_{n=0}^{\infty}$ 
is an adelic sequence of $(X,L)$. 
We call $\{(\calX_n, \overline{\calL_n})\}_{n=0}^{\infty}$ 
the adelic sequence of $(X,L)$ associated with 
$(\calX_0, \overline{\calL_0})$ and $\{f_1, \cdots, f_k\}$. 
\item
Let $(\calX^{'}_0, \overline{\calL^{'}_0})$ be another 
$C^{\infty}$-model of $(X,L)$ 
such that $\calX^{'}_0$ is normal and $\calL^{'}_0$ is nef.  
Let $\{(\calX^{'}_n, \overline{\calL^{'}_n})\}_{n=0}^{\infty}$ 
be the adelic sequence associated with 
$(\calX^{'}_0, \overline{\calL^{'}_0})$ and 
$\{f_1, \cdots, f_k\}$. 
Then $\{(\calX_n, \overline{\calL_n})\}_{n=0}^{\infty}$ 
and $\{(\calX^{'}_n, \overline{\calL^{'}_n})\}_{n=0}^{\infty}$ 
are equivalent.
\end{enumerate}
\end{Theorem}

\Proof
(1) By construction, $\overline{\calL_n}$ is nef for all $n \geq 0$. 
Moreover, over $U$, $\rest{\calX_{n}}{U} = \calX_U$ 
and $\rest{\calL_{n+1}}{U} = \left[f_1^*(\calL_{U})\otimes \cdots 
\otimes f_k^*(\calL_{U})\right]^{\otimes \frac{1}{d}} 
= {\calL}_{U}$. Thus, we have only to check the condition (iii) 
of Definition~\ref{def:adelic:sequence}. 

Let us first estimate the difference between 
$(\calX_1, \overline{\calL_1})$ and $(\calX_2, \overline{\calL_2})$, 
using the difference between $(\calX_0, \overline{\calL_0})$ 
and $(\calX_1, \overline{\calL_1})$. 
Let $\calX_0^{(i)}$ (resp. $\calX_1^{(i)}$) be the normalization of 
$\calX_U \overset{f_i}{\longrightarrow} \calX_U \hookrightarrow \calX_0$ 
(resp.  
$\calX_U \overset{f_i}{\longrightarrow} \calX_U \hookrightarrow \calX_1$).
We denote the induced morphism by $f_0^{(i)}: \calX_0^{(i)} \to \calX_0$  
(resp. $f_1^{(i)}: \calX_1^{(i)} \to \calX_1$). 
By the definition of $\{(\calX_n, \overline{\calL_n})\}_{n=0}^{\infty}$, 
$\calX_1$ (resp. $\calX_2$) is the Zariski closure of 
$\calX_U \overset{\Delta}{\longrightarrow}
\calX_U\times_U\cdots\times_U\calX_U
\hookrightarrow 
\calX_0^{(1)}\times_{O_K}\cdots\times_{O_K}\calX_0^{(k)}$, 
(resp. $\calX_U \overset{\Delta}{\longrightarrow}
\calX_U\times_U\cdots\times_U\calX_U
\hookrightarrow 
\calX_1^{(1)}\times_{O_K}\cdots\times_{O_K}\calX_1^{(k)}$. 
Let $p_{0i}: \calX_0^{(1)}\times_{O_K}\cdots\times_{O_K}\calX_0^{(k)} 
\to \calX_0^{(i)}$ (resp. $p_{1i}: 
\calX_1^{(1)}\times_{O_K}\cdots\times_{O_K}\calX_1^{(k)} 
\to \calX_1^{(i)}$
denote the projection to the $i$-th factor. 

Take a projective arithmetic variety $\calY$ 
such that there exist birational morphisms 
$\rho_0: \calY \to \calX_0$ and $\rho_1: \calY \to \calX_1$
that are the identity maps over $U$. For example, we may take $\calY$ 
as the integral closure of the image of the diagonal map from
$X$ to $\calX_0 \times_{O_K} \calX_1$. 
Let $\calY^{(i)}$ be the normalization of 
\[
\calY_U = \calX_U \overset{f_i}{\longrightarrow} 
\calY_U = \calX_U \hookrightarrow \calY.
\]
We denote the induced morphism by $g^{(i)}: \calY^{(i)} \to \calY$. 
Let $\calZ$ be the Zariski closure of 
\[
\calY_U \overset{\Delta}{\longrightarrow}
\overbrace{\calY_U\times_U\cdots\times_U\calY_U}^{k}
\hookrightarrow 
\calY^{(1)}\times_{O_K}\cdots\times_{O_K}\calY^{(k)}, 
\]
where $\Delta$ is the diagonal map.
Then there are biratinal morphisms 
$\mu_0: \calZ \to \calX_1$ and 
$\mu_1: \calZ \to \calX_2$ that are identity maps over $U$. 
Let $q_i: \calY^{(1)}\times_{O_K}\cdots\times_{O_K}\calY^{(k)} 
\to \calY^{(i)}$ denote the projection to the $i$-th factor. 

\begin{Claim}
\label{claim:1}
\renewcommand{\labelenumi}{(\arabic{enumi})}
\begin{enumerate}
\item
$\displaystyle{\rest{\left(\bigotimes_{i=1}^k 
(\rho_0 \circ g^{(i)}\circ q_i)^*\overline{\calL_0}\right)}{\calZ}
= \mu_0^*\left(\rest{\bigotimes_{i=1}^k(f_0^{(i)} \circ p_{0i})^*
\overline{\calL_0}}{\calX_1}\right)}$. 
\item
$\displaystyle{\rest{\left(\bigotimes_{i=1}^k 
(\rho_1 \circ g^{(i)}\circ q_i)^*\overline{\calL_1}\right)}{\calZ}
= \mu_1^*\left(\rest{\bigotimes_{i=1}^k(f_1^{(i)} \circ p_{1i})^*
\overline{\calL_1}}{\calX_2}\right)}$. 
\end{enumerate}
\end{Claim}

Indeed, for each $i=1, \cdots, k$, we consider two morphisms:
\begin{align*}
& \calZ \overset{\Delta}{\longrightarrow} 
\calY^{(1)}\times_{O_K}\cdots\times_{O_K}\calY^{(k)}
\overset{g^{(i)}\circ q_i}{\longrightarrow} 
\calY \overset{\rho_0}{\longrightarrow} \calX_0 
\qquad\text{and} \\
& \calZ \overset{\mu_0}{\longrightarrow}\calX_1 
\overset{\Delta}{\longrightarrow}
\calX^{(1)}\times_{O_K}\cdots\times_{O_K}\calX^{(k)}
\overset{f_0^{(i)}\circ p_{0i}}{\longrightarrow} \calX_0. 
\end{align*}
These two morphisms coincide with each other over $U$, 
and hence over $O_K$. 
Thus we get the first assertion. 
The second assertion follows by the same argument. 

\begin{Claim}
\label{claim:2}
$\displaystyle{\mu_0^*(\overline{\calL_1}) - \mu_1^*(\overline{\calL_2})
= \frac{1}{d}
\rest{\left(\bigotimes_{i=1}^k 
(\rho_0 \circ g^{(i)}\circ q_i)^*\overline{\calL_0}
- (\rho_1 \circ g^{(i)}\circ q_i)^*\overline{\calL_1}
\right)}{\calZ}}$
\end{Claim}

Indeed, it follows from Claim~\ref{claim:1} that 
\begin{align*}
\rest{\left(\sum_{i=1}^k 
(\rho_0 \circ g^{(i)}\circ q_i)^*\overline{\calL_0}
\right)}{\calZ}
&=
\mu_0^*\left(\rest{\sum_{i=1}^k(f_0^{(i)} \circ p_{0i})^*
\overline{\calL_0}}{\calX_1}\right) \\
&= \mu_0^*(d \overline{\calL_1}) = d \mu_0^*(\overline{\calL_1})
\end{align*}
and similarly 
\[
\rest{\left(\sum_{i=1}^k 
(\rho_1 \circ g^{(i)}\circ q_i)^*\overline{\calL_1}
\right)}{\calZ}
= d \mu_1^*(\overline{\calL_2}). 
\]
Since $\rho_0: \calY \to \calX_0$ and $\rho_1: \calY \to \calX_1$ 
give the same morphism over $U$, there are positive rational numbers  
$c_{01}, \cdots, c_{0r}$ and positive real numbers 
$c_{0\sigma} \;(\sigma:K \to \CC)$ 
such that 
\begin{multline*}
\pi_{\calY}^*\left(\left( - \sum_{\alpha=1}^r c_{0\alpha} 
[\frakP_{\alpha}],
- \sum_{\sigma:K \to \CC} c_{0\sigma} [\sigma]
\right)\right) \lesssim 
\rho_0^*(\overline{\calL_0}) - \rho_1^*(\overline{\calL_{1}}) \\
\lesssim \pi_{\calY}^*\left(\left(  \sum_{\alpha=1}^r c_{0\alpha} 
[\frakP_{\alpha}], 
\sum_{\sigma:K \to \CC} c_{0\sigma} [\sigma]
\right)\right). 
\end{multline*}
Then it follows from Claim~\ref{claim:2} that 
\begin{multline*}
\pi_{\calZ}^*\left(\left(
- \sum_{\alpha=1}^r \frac{k}{d} c_{0\alpha} 
[\frakP_{\alpha}], 
- \sum_{\sigma:K \to \CC} \frac{k}{d} c_{0\sigma} [\sigma]
\right)\right) \lesssim 
\mu_0^*(\overline{\calL_1}) - \mu_1^*(\overline{\calL_{2}}) \\
\lesssim \pi_{\calZ}^*\left(\left(
  \sum_{\alpha=1}^r \frac{k}{d} c_{0\alpha} 
[\frakP_{\alpha}], 
\sum_{\sigma:K \to \CC} \frac{k}{d} c_{0\sigma} [\sigma]
\right)\right). 
\end{multline*}
Inductively we find that there exist 
a projective arithmetic variety $\calW$, and 
birational morphisms $\nu_n: \calW \to \calX_n$ and 
$\nu_n: \calW \to \calX_{n+1}$ such that 
\begin{multline*}
\pi_{\calW}^*\left(\left(
 - \sum_{\alpha=1}^r \frac{k^n}{d^n} c_{0\alpha} 
[\frakP_{\alpha}], 
- \sum_{\sigma:K \to \CC} \frac{k^n}{d^n} c_{0\sigma} [\sigma]
\right)\right) \lesssim 
\mu_n^*(\overline{\calL_n}) - \nu_n^*(\overline{\calL_{n+1}}) \\
\lesssim \pi_{\calW}^*\left(\left(
  \sum_{\alpha=1}^r \frac{k^n}{d^n} c_{0\alpha} 
[\frakP_{\alpha}], 
\sum_{\sigma:K \to \CC} \frac{k^n}{d^n} c_{0\sigma} [\sigma]
\right)\right).
\end{multline*}
Since $\sum_{n=0}^{\infty} \left(\frac{k}{d}\right)^n < \infty$, 
$\{(\calX_n, \overline{\calL_n})\}_{n=0}^{\infty}$ satisfies 
the condition (iii) of Definition~\ref{def:adelic:sequence}. 

\smallskip
\noindent
(2) 
Take a non-empty Zariski open set $\widetilde{U}$ of $\Spec(O_K)$ 
that is a common base of both 
$\{(\calX_n, \overline{\calL_n})\}_{n=0}^{\infty}$
and $\{(\calX^{'}_n, \overline{\calL^{'}_n})\}_{n=0}^{\infty}$. 
Let $\calX^{''}_0$ be a projective arithmetic variety  
such that there exist birational morphisms 
${\widetilde{\mu}}_0: \calX^{''}_0 \to \calX_0$ and 
${\widetilde{\nu}}_0: \calX^{''}_0 \to \calX^{'}_0$
that are the identity maps over $\widetilde{U}$. 
Then there are positive rational numbers 
$\widetilde{c}_{01}, \cdots, \widetilde{c}_{0\widetilde{r}}$ 
and positive real numbers ${\widetilde{c}}_{0\sigma} 
\;(\sigma:K \to \CC)$ such that 
\begin{multline*}
\pi_{\calX^{''}_0}^*\left(\left( - \sum_{\alpha=1}^{\widetilde{r}} 
\widetilde{c}_{0\alpha} 
[\widetilde{\frakP}_{\alpha}], 
- \sum_{\sigma:K \to \CC} {\widetilde{c}}_{0\sigma} [\sigma]
\right)\right) \lesssim 
{\widetilde{\mu}}_0^*(\overline{\calL_0}) - 
{\widetilde{\nu}}_0^*(\overline{\calL^{'}_{0}}) \\
\lesssim \pi_{\calX^{''}_0}^*\left(\left( 
\sum_{\alpha=1}^{\widetilde{r}} 
\widetilde{c}_{0\alpha} 
[\widetilde{\frakP}_{\alpha}], 
\sum_{\sigma:K \to \CC} {\widetilde{c}}_{0\sigma} [\sigma]
\right)\right),  
\end{multline*}
where $\Spec{O_K}\setminus \widetilde{U} = 
\{\widetilde{\frakP}_1, \cdots, \widetilde{\frakP}_{\widetilde{r}}\}$ 
and $\pi_{\calX^{''}_0}: \calX^{''}_0 \to \Spec(O_K)$ is the structure 
morphism. 
Let $\calX^{''(i)}_{0}$ be the normalization of 
$\calX^{''}_{0\widetilde{U}} \overset{f_i}{\longrightarrow} 
\calX^{''}_{0\widetilde{U}} \hookrightarrow \calX^{''}_{0}$. 
Let $\calX^{''}_{1}$ be the Zariski closure of 
$\calX^{''}_{0\widetilde{U}} \overset{\Delta}{\longrightarrow}
{\calX^{''}_{0\widetilde{U}}\times_{\widetilde{U}}
\cdots\times_{\widetilde{U}}\calX^{''}_{0\widetilde{U}}}
\hookrightarrow 
\calX^{''(1)}_{0}\times_{O_K}\cdots
\times_{O_K}\calX^{''(k)}_{0}$.
Then there are biratinal morphisms 
${\widetilde{\mu}}_1: \calX^{''}_{1} \to \calX_1$ and 
${\widetilde{\nu}}_1: \calX^{''}_{1} \to \calX^{'}_1$ 
that are identity maps over $\widetilde{U}$. 
By the same argument as in (1), we have 
\begin{multline*}
\pi_{\calX^{''}_{1}}^*\left(\left( 
- \sum_{\alpha=1}^{\widetilde{r}} \frac{k}{d} 
\widetilde{c}_{0\alpha} [\widetilde{\frakP}_{\alpha}], 
- \sum_{\sigma:K \to \CC} \frac{k}{d} {\widetilde{c}}_{0\sigma} [\sigma]
\right)\right) \lesssim 
{\widetilde{\mu}}_1^*(\overline{\calL_1}) 
- {\widetilde{\nu}}_1^*(\overline{\calL^{'}_{1}}) \\
\lesssim \pi_{\calX^{''}_{1}}^*\left(\left(  
\sum_{\alpha=1}^{\widetilde{r}} \frac{k}{d} \widetilde{c}_{0\alpha} 
[\widetilde{\frakP}_{\alpha}], 
\sum_{\sigma:K \to \CC} \frac{k}{d} {\widetilde{c}}_{0\sigma} [\sigma]
\right)\right). 
\end{multline*}
Inductively we find that there exist 
a projective arithmetic variety $\calX^{''}_n$, 
birational morphisms $\widetilde{\mu}_n: \calX^{''}_n \to \calX_n$ and 
$\widetilde{\nu}_n: \calX^{''}_n \to \calX^{'}_{n}$ such that 
\begin{multline*}
\pi_{\calX^{''}_n}^{*} 
\left(\left( - \sum_{\alpha=1}^{\widetilde{r}} \frac{k^n}{d^n} 
{\widetilde{c}}_{n\alpha} 
[{\widetilde{\frakP}}_{\alpha}], 
- \sum_{\sigma:K \to \CC} \frac{k^n}{d^n} {\widetilde{c}}_{0\sigma} [\sigma]
\right)\right) 
\lesssim
\widetilde{\mu}_n^*(\overline{\calL_n}) - 
\widetilde{\nu}_n^*(\overline{\calL^{'}_{n}})\\ 
\lesssim
\pi_{\calX^{''}_n}^*
\left(\left(  \sum_{\alpha=1}^{\widetilde{r}} \frac{k^n}{d^n} 
{\widetilde{c}}_{0\alpha} 
[{\widetilde{\frakP}}_{\alpha}], 
\sum_{\sigma:K \to \CC}  \frac{k^n}{d^n} {\widetilde{c}}_{0\sigma} [\sigma]
\right)\right). 
\end{multline*}
Since $\sum_{n=0}^{\infty} \frac{k^n}{d^n} < \infty$, 
we get the assertion. 
\QED

\subsection{Adelic intersection and heights of subvarieties}
As in \S\ref{subsec:adelic:sequence:from:morphisms}, 
let $(X; f_1, \cdots, f_k)$ be a dynamical system of $k$ morphisms 
over $K$ associated with $L$ of degree $d >k $, 
and we assume that $X$ is normal, $f_1, \cdots, f_k$ are surjective, 
and $L$ is an ample $\QQ$-line bundle. 
Fix a $C^{\infty}$-model $(\calX, \overline{\calL})$ 
of $(X,L)$. 

Let $Y \subset X_{\overline{K}}$ be a subvariety. 
Take a finite extension field $K'$ of $K$ over which  
$Y$ is defined. Let $\calY$ denote the Zariski closure 
of $Y$ in $\calX\otimes_{O_K}O_{K'}$. We define the height of 
$Y$ with respect to $(\calX, \overline{\calL})$ to be 
\[
{h}_{(\calX, \overline{\calL})}(Y)
= 
\frac{\adeg\left(\achern_1
\left(\rest{\overline{\calL}\otimes_{O_K}O_{K'}}{\calY}
\right)^{\dim Y +1}\right)}%
{[K':\QQ](\dim Y +1)\deg(\rest{L\otimes_K K'}{Y}^{\cdot \dim Y})}.
\] 
In particular, if $Y$ is a closed point of $X_{\overline{K}}$, 
then 
\[
{h}_{(\calX, \overline{\calL})}(Y)
= 
\frac{\adeg\left(\rest{\overline{\calL}\otimes_{O_K}O_{K'}}
{\calY}\right)}{[K':\QQ]}.
\]
Then we have the following theorem.

\begin{Theorem}
\label{thm:canonical:height:of:subvarieties} 
Let $(X; f_1, \cdots, f_k)$ be a dynamical system of $k$ morphisms 
over $K$ associated with $L$ of degree $d >k $.  
We assume that $X$ is normal, $f_1, \cdots, f_k$ are surjective, 
and $L$ is an ample $\QQ$-line bundle. 
\renewcommand{\labelenumi}{(\arabic{enumi})}
\begin{enumerate}
\item 
Let $(\calX_0, \overline{\calL_0})$ be a $C^{\infty}$-model of 
$(X,L)$ such that $\calX_0$ is normal and $\overline{\calL_0}$ is nef. 
Let $\{(\calX_n, \overline{\calL_n})\}_{n=0}^{\infty}$ be 
the adelic sequence of $(X,L)$ associated with 
$(\calX_0, \overline{\calL_0})$ and $\calF=\{f_1, \cdots, f_k\}$. 
Then, for any subvariety $Y \subset X_{\overline{K}}$, 
\[
\widehat{h}_{L, \calF}(Y) 
:= \lim_{n\to\infty}{h}_{(\calX_n, \overline{\calL_n})}(Y)
\]
converges. The limit $\widehat{h}_{L, \calF}(Y)$
is independent of the choice of $C^{\infty}$-models 
$(\calX_0, \overline{\calL_0})$ of $(X, L)$, 
and $\widehat{h}_{L, \calF_1}(Y) \geq 0$ 
for all $Y \subset X_{\overline{K}}$.
We call this limit the {\em canonical height} of $Y$ 
associated with $L$ and $\calF$.
\item 
Let $(\calX, \overline{\calL})$ be a $C^{\infty}$-model of $(X,L)$. 
Then there exists a constant $C$ such that 
\[
\left\vert 
\widehat{h}_{L, \calF}(Y) - {h}_{(\calX, \overline{\calL})}(Y)
\right\vert \leq C
\]
for any subvariety $Y \subset X_{\overline{K}}$. 
\item
When $Y$ is a closed point of $X_{\overline{K}}$, 
$\widehat{h}_{L, \calF}(Y)$ coincides with the canonical height 
in Theorem~\ref{thm:main}.
\end{enumerate}
\end{Theorem}

\Proof
(1) Let $\calY_n$ denote the Zariski closure 
of $Y$ in $\calX_n\otimes_{O_K}O_{K'}$. 
By Lemma~\ref{lemma:adelic:sequence:associated}, 
$\{(\calY_n, \rest{\overline{\calL_n}\otimes_{O_K}O_{K'}}{\calY_n})
\}_{n=0}^{\infty}$ is an adelic sequence of $(Y, \rest{L}{Y})$. 
Then it follows from Theorem~\ref{thm:adelic:intersection} 
that ${h}_{(\calX_n, \overline{\calL_n})}(Y)$ converges as $n$ 
goes to $\infty$. 
It follows from Theorem~\ref{thm:adelic:intersection} and 
Theorem~\ref{thm:adelic:sequence:from:morphisms}(2) 
that this limit is independent of the choice of 
$(\calX_0, \overline{\calL_0})$. 
Moreover, 
since $\overline{\calL_n}$ is nef, 
${h}_{(\calX_n, \overline{\calL_n})}(Y) \geq 0$ by 
Proposition~\ref{prop:nef:Q:effective}. 
Thus  $\widehat{h}_{L}(Y) = 
\lim_{n\to\infty}{h}_{(\calX_n, \overline{\calL_n})}(Y) \geq 0$.

\smallskip
(2)
Take a $C^{\infty}$-model $(\calX', \overline{\calL'})$ of 
$(X,L)$ such that $\calX'$ is normal and $\overline{\calL'}$ is nef. 
By \cite{BGS},~Proposition~3.2.2, 
there exists a constant $C_1$ such that 
\[
\left\vert
{h}_{(\calX, \overline{\calL})}(Y) 
- {h}_{(\calX', \overline{\calL'})}(Y)
\right\vert \leq C_1
\]
for any subvariety $Y \subset X_{\overline{K}}$.
Thus, to prove (2), we may assume that 
$(\calX, \overline{\calL})$ is nef. 

We take birational morphisms $\mu_n: \calX_n^{'} \to \calX_n$ and 
$\nu_n: \calX_n^{'} \to \calX_{n+1}$ as in 
Definition~\ref{def:adelic:sequence}. 
Let $\calY'$ be the Zariski closure of $Y$ in $\calX_n^{'}$. 
Then by a similar argument as in the proof of 
Theorem~\ref{thm:adelic:intersection}, we have 
\begin{multline*}
\left\vert
\achern_1(\rest{\mu_n^*(\overline{\calL_n})
\otimes_{O_K} O_{K'}}{\calY'})^{\dim Y +1}
- \achern_1(\rest{\nu_n^*(\overline{\calL_{n+1}})
\otimes_{O_K} O_{K'}}{\calY'})^{\dim Y+1}
\right\vert  \\
\leq
\frac{k^n}{d^n} [K':K] 
\left(\sum_{\alpha=1}^r c_{0\alpha} 
\log\#(O_K/\frakP_{\alpha})
+ \frac{1}{2} \sum_{\sigma: K \to \CC} c_{0\sigma}\right). 
\end{multline*}
Since by the projection formula $\achern_1(\rest{\overline{\calL_n}
\otimes_{O_K} O_{K'}}{\calY_n})^{\dim Y +1} = 
\achern_1(\rest{\mu_n^*(\overline{\calL_n})
\otimes_{O_K} O_{K'}}{\calY'})^{\dim Y +1}$ and 
$\achern_1(\rest{\overline{\calL_{n+1}}
\otimes_{O_K} O_{K'}}{\calY_{n+1}})^{\dim Y +1} = 
\achern_1(\rest{\nu_n^*(\overline{\calL_{n+1}})
\otimes_{O_K} O_{K'}}{\calY'})^{\dim Y +1}$, 
we have 
\[
\left\vert
h_{(\calX_n, \calL_n)}(Y) -  
h_{(\calX_{n+1}, \calL_{n+1})}(Y)
\right\vert
\leq
\frac{ {k^n} 
\left( \sum_{\alpha=1}^r c_{0\alpha} 
\log\#(O_K/\frakP_{\alpha})
+ \frac{1}{2} \sum_{\sigma \to \CC} c_{0\sigma}\right)}%
{d^n [K:\QQ] (\dim Y+1)}. 
\]
Putting $C = 
\frac{ d 
\left( 
\sum_{\alpha=1}^r c_{0\alpha} 
\log\#(O_K/\frakP_{\alpha})
+ \frac{1}{2} \sum_{\sigma \to \CC} c_{0\sigma}\right)}{(d-k)[K:\QQ]}$, 
we obtain the desired estimate. 

\smallskip
(3)
Let us check that if $y \in X(\overline{K})$ is a closed point, 
$\widehat{h}_{L, \calF_1}(y)$ defined via adelic intersection 
satisfies (i) and (ii) of Theorem~\ref{thm:main}. 
Note that $h_{(\calX_0, \overline{\calL_0})} = h_L + O(1)$. 
Then it follows from (2) that  
$\widehat{h}_{L, \calF} = h_L + O(1)$. To show (ii), 
let $K'$ be a finite extension field of $K$ over which
$y$ is defined. Let $\Delta_y^{(n)}$ 
denote the Zariski closure of $y$ in $\calX_n \times_{\Spec(O_K)} 
\Spec(O_{K'})$. Put 
\[
b_n(y) : = 
\frac{\adeg\left(\rest{\overline{\calL_{n}}
\otimes_{O_K} O_{K'}}{\Delta_y^{(n)}}\right)}%
{[K':\QQ]}. 
\]
Then $b_n(y)$ is independent of the choice of $K'$, 
and by definition 
$\widehat{h}_{L, \calF}(y) = \lim_{n \to \infty} b_n(y)$. 

\begin{Claim}
\label{claim:4}
$\sum_{i=1}^{k} b_n(f_i(y)) 
= d \, b_{n+1}(y)$. 
\end{Claim}

Let us show the assertion for $n=0$. The assertion for $n \geq 1$ 
can be shown by the same argument. We follow the notation 
in the beginning of \S\ref{subsec:adelic:sequence:from:morphisms}. 
By the projection formula, we have 
\begin{align*}
  & \sum_{i=1}^{k} b_0(f_i(y)) - d \, b_{1}(y) \\
= & \sum_{i=1}^k \frac{\adeg\left(\rest{\overline{\calL_{0}}
\otimes_{O_K} O_{K'}}{\Delta_{f_i(y)}^{(0)}}\right)}{[K':\QQ]}
- d \frac{\adeg\left(\rest{\overline{\calL_{1}}
\otimes_{O_K} O_{K'}}{\Delta_y^{(1)}}\right)}{[K':\QQ]} \\
= & \sum_{i=1}^k \frac{\adeg\left(\rest{(f^{(i)}\circ p_i)^* 
(\overline{\calL_{0}}\otimes_{O_K} O_{K'})}{\Delta_{y}^{(1)}}\right)}{[K':\QQ]}
- d \frac{\adeg\left(\rest{\overline{\calL_{1}}
\otimes_{O_K} O_{K'}}{\Delta_y^{(1)}}\right)}{[K':\QQ]} \\
= & \frac{\adeg\left( \sum_{i=1}^ k \rest{(f^{(i)}\circ p_i)^* (\overline{\calL_{0}}\otimes_{O_K} O_{K'}) - d (\overline{\calL_{1}}
\otimes_{O_K} O_{K'})}{\Delta_{y}^{(1)}}\right)}{[K':\QQ]}
= 0 
\end{align*}

By letting $n$ to $\infty$ in Claim~\ref{claim:4}, 
we get 
$\sum_{i=1}^{k} \widehat{h}_{L, \calF} (f_i(y)) 
= {d} \widehat{h}_{L, \calF}(y)$.  
\QED

\medskip
\section{Admissible metrics and invariant currents}
\label{sec:metrics:currents}
In this section, we introduce an invariant current for a dynamical system of
$k$ morphisms over $\CC$ associated with a line bundle $L$ of degree $d >
k$.
If $X=\PP^N$, $f:\PP^N \to \PP^N$ is a morphism of degree $d \geq 2$, and
$L= \OO_{\PP^N}(1)$, then this invariant current coincides with the Green
current of $f$ introduced by Hubbard--Papadopol \cite{HP}. 
When $(X;\sigma_1, \sigma_2)$ is Wheler's K3 surface, we will relate the
invariant current with the currents $T^+$ and $T^-$ corresponding
respectively to $\sigma_2\circ\sigma_1$ and $\sigma_1\circ\sigma_2$
introduced by Cantat \cite{Ca}. 
Finally we pose a question regarding distribution of small points, 
and show that this question is true for Latt\'es examples, which are 
certain endomorphisms of $\PP^n$. 

\subsection{Admissible metrics}
\label{subsec:admissible:metrics}
Let $(X; f_1, \cdots, f_k)$ be a
dynamical system of $k$ morphisms over $\CC$ associated with a line bundle $L$
of degree $d > k$. We fix an isomorphism
$\varphi: L^{\otimes d} \overset{\sim}{\rightarrow}
f_1^*L\otimes\cdots\otimes f_k^*L$.
The following theorem is a generalization of \cite{Zh},~Theorem(2.2) of one
morphism to $k$ morphisms over $\CC$.

\begin{Theorem}
\label{thm:admissible:metrics}
\begin{enumerate}
\item[(1)]
Let $\Vert\cdot\Vert_0$ be a continuous metric on $L$.
We inductively define a metric $\Vert\cdot\Vert_n$ on $L$ by
\[
\Vert\cdot\Vert_{n}^d =
\varphi^*(f_1^*\Vert\cdot\Vert_{n-1} \cdots f_k^*\Vert\cdot\Vert_{n-1}).
\]
Then $\Vert\cdot\Vert_n$ converges uniformly
to a metric, which we denote by $\Vert\cdot\Vert_{\infty}$.
\item[(2)]
$\Vert\cdot\Vert_{\infty}$ is the unique continuous metric on $L$ with
\[
\Vert\cdot\Vert_{\infty}^d =
\varphi^*(f_1^*\Vert\cdot\Vert_{\infty} \cdots
f_k^*\Vert\cdot\Vert_{\infty}).
\]
\item[(3)]
If $\varphi$ is replaced with $c\varphi$ 
\textup{(}$c\neq 0 \in \CC$\textup{)},
then $\Vert\cdot\Vert_{\infty}$ is replaced with
$|c|^{\frac{1}{d-k}} \Vert\cdot\Vert_{\infty}$.
\end{enumerate}
We call the metric $\Vert\cdot\Vert_{\infty}$ the 
{\em admissible metric} for $L$. 
\end{Theorem}

\Proof
Set $H = \frac{\Vert\cdot\Vert_{1}}{\Vert\cdot\Vert_{0}}$. Then $H$ is a
continuous function on $X(\CC)$. We have
\begin{align*}
\frac{\Vert\cdot\Vert_{n+1}}{\Vert\cdot\Vert_{n}}(x)
& =
\left(
\frac{\varphi^*(f_1^*\Vert\cdot\Vert_{n} \cdots f_k^*\Vert\cdot\Vert_{n})}
{\varphi^*(f_1^*\Vert\cdot\Vert_{n-1} \cdots f_k^*\Vert\cdot\Vert_{n-1})}
\right)^{\frac{1}{d}}(x) \\
& = 
\left(
\frac{\Vert\cdot\Vert_{n}}{\Vert\cdot\Vert_{n-1}}(f_1(x))
\cdots
\frac{\Vert\cdot\Vert_{n}}{\Vert\cdot\Vert_{n-1}}(f_k(x))
\right)^{\frac{1}{d}}\\
& = \cdots =
\left(
\prod_{f \in \calF_n}
\frac{\Vert\cdot\Vert_{1}}{\Vert\cdot\Vert_{0}}(f(x))
\right)^{\frac{1}{d^n}}
=
\left(
\prod_{f \in \calF_n} H \circ f
\right)^{\frac{1}{d^n}}(x).
\end{align*}
Thus we get
$\Vert\cdot\Vert_{n}
=
\Vert\cdot\Vert_{0} \cdot
\prod_{l=0}^{n-1}
\left(
\prod_{f \in \calF_l} H \circ f
\right)^{\frac{1}{d^l}}$.
We set  $H_n = \prod_{l=0}^{n-1} \prod_{f \in \calF_l}
\left(H \circ f\right)^{\frac{1}{d^l}}$. 

\begin{Claim}
The function $H_n$ converges uniformly as $n \to \infty$.
\end{Claim}

For positive integers $m$ and $n$, we have 
\[
H_{n+m} - H_n 
=
\left(
\prod_{l=0}^{n-1} \prod_{f \in \calF_l} 
\left(H \circ f\right)^{\frac{1}{d^l}}
\right)
\left(
\prod_{l=n}^{m-1} \prod_{f \in \calF_l} 
\left(H \circ f\right)^{\frac{1}{d^l}}
-1 \right). 
\]
We set $M =\max\{1, \sup_{x \in X(\CC)} H(x)\}$ 
and $m = \min\{1, \inf_{x \in X(\CC)} H(x)\}$. Then 
$0 < m \leq M$. It follows from  
$\prod_{f \in \calF_l}\left(H \circ f\right)^{\frac{1}{d^l}} \leq M^{(\frac{k}{d})^{l}}$ 
that 
\[
\prod_{l=0}^{n-1} \prod_{f \in \calF_l} 
\left(H \circ f\right)^{\frac{1}{d^l}}
\leq M^{\sum_{l=0}^{n-1}(\frac{k}{d})^{l}}
\leq M^{\sum_{l=0}^{\infty}(\frac{k}{d})^{l}}
= M^{\frac{d}{d-k}}. 
\] 
We also have 
$\prod_{l=n}^{m-1} \prod_{f \in \calF_l} 
\left(H \circ f\right)^{\frac{1}{d^l}}
 \leq M^{\sum_{l=n}^{m-1}(\frac{k}{d})^{l}}
 \leq M^{(\frac{k}{d})^{n}\frac{d}{d-k}}$ and 
$\prod_{l=n}^{m-1} \prod_{f \in \calF_l} 
\left(H \circ f\right)^{\frac{1}{d^l}}
 \geq m^{\sum_{l=n}^{m-1}(\frac{k}{d})^{l}}
 \geq m^{(\frac{k}{d})^{n}\frac{d}{d-k}}$. 
Then we have 
\[
\Vert H_{n+m} - H_n \Vert_{\sup} \leq 
M^{\frac{d}{d-k}} \max\left\{
M^{(\frac{k}{d})^{n}\frac{d}{d-k}} - 1, 
1 - m^{(\frac{k}{d})^{n}\frac{d}{d-k}}
\right\} 
\rightarrow 
0 \quad
(n, m \to\infty). 
\]
This shows $H_n$ converges uniformly. 
Let $H_{\infty}$ denote 
the limit of $H_n$. 

Since 
\[
\inf_{x \in X(\CC)} H_{\infty}(x)
= \inf_{x \in X(\CC)} \prod_{l=0}^{\infty} \prod_{f \in \calF_l} 
\left(H(f(x))\right)^{\frac{1}{d^l}}
\geq m^{\sum_{l=0}^{\infty}(\frac{k}{d})^l}
\geq m^{\frac{d}{d-k}}, 
\]
we find that $H_{\infty}$ is a positive function.  
Thus $\Vert\cdot\Vert_{n} = \Vert\cdot\Vert_{0} H_n$ 
converges uniformly to a metric 
$\Vert\cdot\Vert_{\infty}:= \Vert\cdot\Vert_{0} H_{\infty}$. 
Thus we have shown (1). 

Next we show (2). 
Letting $n \to\infty$ in
$\Vert\cdot\Vert_{n+1}^d =
\varphi^*(f_1^*\Vert\cdot\Vert_{n} \cdots f_k^*\Vert\cdot\Vert_{n})$,
we get
$\Vert\cdot\Vert_{\infty}^d =
\varphi^*(f_1^*\Vert\cdot\Vert_{\infty} \cdots
f_k^*\Vert\cdot\Vert_{\infty})$. 
To show the uniqueness of $\Vert\cdot\Vert_{\infty}$, 
suppose continuous metrics 
$\Vert\cdot\Vert_{\infty}$ and $\Vert\cdot\Vert_{\infty}^{'}$ 
both satisfy the equality in (2). We set 
$A = \frac{\Vert\cdot\Vert_{\infty}^{'}}{\Vert\cdot\Vert_{\infty}}$. 
Then $A$ is a positive continuous function on $X(\CC)$, 
and satisfies $A^d = \prod_{i=1}^k A\circ f_i$. 
Then it follows from 
$\sup_{x \in X(\CC)} A(x) = 
\sqrt[d]{\sup_{x \in X(\CC)} \prod_{i=1}^k A(f_i(x))}
\leq \left(\sup_{x \in X(\CC)} A(x) \right)^{\frac{k}{d}}$
and $\frac{k}{d} <1$ that $\sup_{x \in X(\CC)} A(x) \leq 1$. 
It follows from 
$\inf_{x \in X(\CC)} A(x) = 
\sqrt[d]{\inf_{x \in X(\CC)} \prod_{i=1}^k A(f_i(x))}
\geq \left(\inf_{x \in X(\CC)} A(x) \right)^{\frac{k}{d}}$
and $\frac{k}{d} <1$ that $\inf_{x \in X(\CC)} A(x) \geq 1$.
Hence $A$ is identically $1$. We have shown (2). 

To show (3), set
$\widetilde{\varphi}=c \varphi$.
Let $\widetilde{\Vert\cdot\Vert_{\infty}}$ denote
the metric with respect to $c \varphi$.
Set $\alpha \Vert\cdot\Vert_{\infty} =
\widetilde{\Vert\cdot\Vert_{\infty}}$.
Since 
$\widetilde{\varphi}(f_1^*\widetilde{\Vert\cdot\Vert_{\infty}} \cdots
f_k^*\widetilde{\Vert\cdot\Vert_{\infty}}) =
\widetilde{\Vert\cdot\Vert_{\infty}}^{d}$, 
we have $|c| \alpha^k = \alpha^d$. Thus $\alpha= |c|^{\frac{1}{d-k}}$.
\QED

\subsection{Invariant currents}
\label{subsec:invariant:currents}
Let $M$ be a complex manifold. A closed $(1,1)$-current $S$ on $M$ is said
to admit a locally continuous potential if, for every $m\in M$, there exists
an open neighborhood $U$ of $m$ and a continuous function $u$ on $U$ such
that
$S = dd^c u$.

Let $\pi: N \to M$ be a morphism of complex manifolds. If a closed
$(1,1)$-current $S$ on $M$ admits a locally continuous potential, we can
define its pull-back $\pi^*S$ by $\rest{\pi^*S}{\pi^{-1}(U)} = dd^c \pi^*u$.
Since $\pi^*u$ is continuous and unique up to a pluriharmonic function,
this formula is well-defined and defines a global current $\pi^*S$ on $N$.
(cf. \cite{HP},~p330.)

We denote by $A^{p,p}(M)$ the space of $C^{\infty}$ $(p,p)$-forms on $X$,
and by $D^{p,p}(M)$ the space of $(p,p)$-currents on $X$.
For $\eta \in A^{p,p}(M)$, we denote by $[\eta] \in D^{p,p}(X)$ the current
corresponding to $\eta$.

The next theorem is the main theorem of this section.

\begin{Theorem}
\label{thm:invariant:currents}
Let $(X; f_1, \cdots, f_k)$ be a dynamical system of $k$ morphisms over
$\CC$ associated with a line bundle $L$ of degree $d > k$. We assume that
$X$ is smooth and $L$ is ample.
\begin{enumerate}
\item[(1)]
Let $\Vert\cdot\Vert_{\infty}$ be the admissible metric of $L$.
We set $T = \chern_1(L, \Vert\cdot\Vert_{\infty}) :=
dd^c [-\log \Vert s \Vert_{\infty}^2] + \delta_{\zero(s)}
\in D^{1,1}(X(\CC))$, where $s$ is a non-zero rational section of $L$. 
\textup{(}This
formula does not depend on $s$ and defines a closed $(1,1)$-current on
$X$.\textup{)}
Then $T$ is positive.
\item[(2)]
Since $T$ admits a locally continuous potential by \textup{(1)},
we have the pull-back $f_i^*T \in D^{1,1}(X(\CC))$. Then we have
\[
f_1^*T + \cdots + f_k^*T = dT.
\]
\item[(3)]
Let $\eta_0 \in A^{1,1}(X(\CC))$ be any closed $C^{\infty}$ $(1,1)$-form
whose cohomology class coincides with $\chern_1(L)$. We inductively define
$\eta_n \in A^{1,1}(X(\CC))$ by
\[
\eta_{n+1} = \frac{1}{d} (f_1^*\eta_n + \cdots + f_k^*\eta_n).
\]
Then, $[\eta_n]$ converges to $T$ as currents.
\end{enumerate}
\end{Theorem}

\Proof
(1)
Since $X$ is a smooth projective variety and $L$ is an ample line bundle,
there exists a $C^{\infty}$ metric $\Vert\cdot\Vert_0$ on $L$ such that
its Chern form $\omega_0 := \chern_1(L, \Vert\cdot\Vert_{\infty})
\in A^{1,1}(X)$ is everywhere positive.

By fixing $\varphi: L^{\otimes d} \overset{\sim}{\rightarrow}
f_1^*L_1\otimes\cdots\otimes f_k^*L_k$,
as in the proof of Theorem~\ref{thm:admissible:metrics},
we inductively define a metric $\Vert\cdot\Vert_n$ on $L$ by
\[
\Vert\cdot\Vert_{n+1} =
\left\{
\varphi^*(f_1^*\Vert\cdot\Vert_{n} \cdots f_k^*\Vert\cdot\Vert_{n})
\right\}^{\frac{1}{d}}.
\]
As in Theorem~\ref{thm:admissible:metrics}, 
let $H = \frac{\Vert\cdot\Vert_{1}}{\Vert\cdot\Vert_{0}}$. 
Then $H$ is a $C^{\infty}$ positive function on $X(\CC)$. 
By replacing $\varphi$ with $c \varphi$ with $c >0$
if necessary, we may assume that $H \geq 1$. (Indeed we may take $c =
(\min_{x \in X(\CC)} H(x))^{\frac{1}{d}}$.)

We define a $C^{\infty}$ $(1,1)$-form $\omega_n$ by 
$\omega_n = \frac{1}{d^n} \sum_{f \in \calF_n}
f^*\omega_0$. 

For $x \in X(\CC)$, we take a non-zero rational section $s$ of $L$ such
that $s(x) \neq 0$. We set $U = X\setminus \Supp(\zero(s))$.
We define a $C^{\infty}$ map $G_{n,s}: U \to \RR$ by
$G_{n,s} = -\log\Vert s\Vert_n^2$.

\begin{Claim}
\label{claim:invariant:currents:1}
$G_{n,s}$ is non-increasing with respect to $n$.
Moreover, $dd^c G_{n,s} =\rest{\omega_n}{U}$, 
and thus $dd^c G_{n,s}$ is everywhere positive.
\end{Claim}

Indeed, by the proof of Theorem~\ref{thm:admissible:metrics}, we have
$\Vert\cdot\Vert_{n+1}
=
\Vert\cdot\Vert_{n} \cdot
\left(\prod_{f \in \calF_n} H \circ f \right)^{\frac{1}{d^n}}$.
Thus we find
\[
-\log \Vert s\Vert_{n+1}^2
=
-\log \Vert s \Vert_{n}^2 +
\frac{1}{d^n} \sum_{f \in \calF_n} - 2 \log H \circ f.
\]
Since $H \geq 1$,
$\frac{1}{d^n} \sum_{f \in \calF_n} - 2 \log H \circ f \leq 0$.
Thus $G_{n,s}$ is non-increasing.

On the other hand, we have
$dd^c G_{n,s} = \rest{\chern_1(L, \Vert\cdot\Vert_{n})}{U}$.
Since 
\begin{align*}
\chern_1(L, \Vert\cdot\Vert_{n}) 
& = \frac{1}{d} \left\{
f_1^* \chern_1(L, \Vert\cdot\Vert_{n-1}) + \cdots + 
f_k^* \chern_1(L, \Vert\cdot\Vert_{n-1}) \right\}\\
& = \cdots = \frac{1}{d^n} \sum_{f \in \calF_n}
f^* \chern_1(L, \Vert\cdot\Vert_{0})  = \omega_n,
\end{align*}
we have $dd^c G_{n,s} =
\rest{\omega_n}{U}$, which is everywhere positive.

Since $G_{n,s}$ is a $C^{\infty}$-function with 
$d d^c G_{n,s} \geq 0$ by Claim~\ref{claim:invariant:currents:1}, 
$G_{n,s}$ is a psh (plurisubharmonic) function.
Since $G_{n,s}$ is non-increasing and bounded from below,
$G_{\infty, s} := \lim_{n \to\infty} G_{n,s}$ is
also a psh function. By Theorem~\ref{thm:admissible:metrics},
we have $G_{\infty, s} = -\log \Vert s \Vert_{\infty}$.
Thus $G_{\infty, s}$ is a continuous psh function, and so
$d d^c [-\log \Vert s \Vert_{\infty}^2 ]$ is a positive current on $U$.

Since positivity of a current is a local condition
and we see that $\chern_1(L, \Vert\cdot\Vert_{\infty})$
does not depend on $s$, we find that
the closed $(1,1)$-current $T = \chern_1(L, \Vert\cdot\Vert_{\infty})$
admits a locally continuous potential and is positive.
We note $\lim_{n\to\infty}[\omega_n] = T$ as currents. 

(2)
Since $T$ admits a locally continuous potential, the currents
$f_i^*T$ are defined. Locally, we have
\begin{align*}
f_1^*T + \cdots + f_1^*T
& = d d^c [f_1^*(-\log \Vert s \Vert_{\infty}^2) + \cdots +
f_k^*(-\log \Vert s \Vert_{\infty}^2)] \\
& = d d^c [-\log \Vert s \Vert_{\infty}^{2d}]
= dT.
\end{align*}
Thus we get (2).

(3)
Let $\eta_0 \in A^{1,1}(X)$ be a closed $(1,1)$-form
whose cohomology class coincides with $\chern_1(L)$.
Since the cohomology class of $\omega_0$ is the same as 
that of $\eta_0$, by $d d^c$-lemma, 
there exists $u \in C^{\infty}(X(\CC))$ with
$\eta_0 - \omega_0 = d d^c u$.
Then we have
\begin{align*}
\eta_n 
& = \left(\frac{1}{d}\right)^n \sum_{f \in \calF_n} f^*(\eta_0) \\
& = \left(\frac{1}{d}\right)^n \sum_{f \in \calF_n} f^*(\omega_0) 
+ d d^c \left(\left(\frac{1}{d}\right)^n \sum_{f \in \calF_n} f^*(u)\right)
= \omega_n 
+ d d^c \left(\left(\frac{1}{d}\right)^n \sum_{f \in \calF_n} f^*(u)\right).
\end{align*}
Here $u$ is a bounded function, and 
$\left\Vert
\left(\frac{1}{d}\right)^n \sum_{f \in \calF_n} f^*(u)
\right\Vert_{\sup}
\leq \left(\frac{k}{d}\right)^n \Vert u\Vert_{\sup}
\rightarrow 0$ as $n \to\infty$. 
Thus, as currents, $\left[
\left(\frac{1}{d}\right)^n \sum_{f \in \calF_n} f^*(u)
\right]$ goes to zero as $n \to \infty$. 
Since $[\omega_n] \to T$ as currents, we get 
$\lim_{n\to\infty} [\eta_n] = T$. 
\QED

\subsection{Examples}
\subsubsection{Projective spaces}
We remark that, for the projective space, the invariant current in
Theorem~\ref{thm:invariant:currents} coincides with the Green current.

\begin{Proposition}
\label{prop:current:projective:space}
Let $f: \PP^N_{\CC} \to \PP^N_{\CC}$ be a morphism over $\CC$
of degree $d \geq 2$, and let $T_G$ denote the Green current of $f$.
We regard $(\PP^N_{\CC}, f)$ as a dynamical system associated with
$\OO_{\PP^N}(1)$ of degree $d$, and we denote by $T$ the invariant
current in Theorem~\ref{thm:invariant:currents}. Then $T = T_G$,
\end{Proposition}

\Proof
By definition, $T_G$ is the limit of $\frac{1}{d^n}(f^n)^* \omega_{FS}$
as $d \to \infty$, where $\omega_{FS}$ is the Fubini-Study metric.
On the other hand, by Theorem~\ref{thm:invariant:currents}(3),
$\frac{1}{d^n}(f^n)^* \omega_{FS}$ converges to $T$. Thus $T = T_G$.
\QED

\begin{Remark}
We remark that our assumption $f$ being a morphism is rather
strong, and that, to define the Green current, $f$ suffices to be an
algebraically stable rational map (cf. \cite{Sib},~Th\'eor\`eme~1.6.1).
\end{Remark}

\subsubsection{Wheler's K3 surfaces}
Let $(S; \sigma_1, \sigma_2)$ be Wheler's K3 surface as in
\S\ref{subsubsec:K3:I}.

\begin{Lemma}
\label{lem:current:K3}
Let $(\sigma_2\circ\sigma_1)^*:
H^{1,1}(S, \RR) \to H^{1,1}(S, \RR)$ be the $\RR$-linear map
induced by $\sigma_2\circ\sigma_1: S \to S$,
and let $\lambda = \max\{ |t| \mid
 \det(t I - (\sigma_2\circ\sigma_1)^*) =0 \}$ be the spectral
 radius of $\sigma_2\circ\sigma_1$. Then $\lambda=7+4\sqrt{3}$.
\end{Lemma}

\Proof
With the notation in \S\ref{subsubsec:K3:I},
$(\sigma_2\circ\sigma_1)^*(E^{+}) = (7+4\sqrt{3}) E^{+}$.
Thus $\lambda \geq 7+4\sqrt{3}$.
On the other hand,
it follows from \cite{Mc},~Theorem~3.2 that 
$(\sigma_2\circ\sigma_1)^*$ has at most one 
eigenvalue $\lambda$ with $|\lambda| > 1$. 
Thus $\lambda = 7+4\sqrt{3}$. 
\QED

It is known that the topological entropy of
$\sigma_2\circ\sigma_1$, 
denoted by $h_{\topo}(\sigma_2\circ\sigma_1)$,  
is equal to $\log\lambda$ 
(cf. \cite{Ca},~Th\'eor\`eme~2.1). 
Cantat \cite{Ca} proved the following theorem about K3 automorphisms 
with strictly positive topological entropy. 

\begin{Theorem}[\cite{Ca}, \S2, cf. \cite{Mc}, Theorem~11.3]
\label{thm:Cantat}
Let $X$ be a K3 surface with an automorphism $\phi$ whose 
topological entropy $h_{\topo}(\phi) = \log\lambda$ is strictly positive. 
Then there exists a closed positive $(1,1)$-current $T^{+}$, 
unique up to scale, such that $f^{*}(T^{+}) 
= \lambda T^{+}$. 
Moreover, for any other closed positive $(1,1)$-current 
$\omega$, there is a constant $c^{+}$ such that 
\[
\frac{1}{\lambda^n} (f^*)^n (\omega) 
\to c^{+} T^{+} \quad (n \to \infty)
\]
as $(1,1)$-currents.  
\end{Theorem}

Considering $\phi^{-1}$, we have a closed positive $(1,1)$-current 
$T^{-}$. The next proposition relates $T^{+}$ and 
$T^{-}$ with the invariant current $T$ in 
Theorem~\ref{thm:invariant:currents} for $(S; \sigma_1, \sigma_2)$. 

\begin{Proposition}
\label{prop:current:Wheler:K3}
Let the notation be as in \S\ref{subsubsec:K3:I}. 
Let $\lambda = 7+4\sqrt{3}$. 
Let $(S; \sigma_1, \sigma_2)$ be Wheler's K3 surface over $\CC$, 
which we regard
as a dynamical system of two morphisms associated
with $L$ of degree $4$. Let $T$ be the invariant current
in Theorem~\ref{thm:invariant:currents}.
Let $\eta$ be a positive $C^{\infty}$ $(1,1)$-form on $S$ 
whose cohomology class is equal to $\chern_1(L)$. We set 
\[
T^{+} = \lim_{n\to\infty} \frac{1}{\lambda^n} 
(\sigma_2\circ\sigma_1)^{*n} (\eta), \qquad 
T^{-} = \lim_{n\to\infty} \frac{1}{\lambda^n} 
(\sigma_2\circ\sigma_1)^{-1*n} (\eta),
\]
whose convergence is assured by Theorem~\ref{thm:Cantat}. 
Then $T = T^+ + T^-$.
\end{Proposition}

\Proof
We set
\begin{align*}
S_n & = \frac{1}{\lambda^n} 
(\sigma_2\circ\sigma_1)^{*n} (\eta)
+ \frac{1}{\lambda^n} 
(\sigma_2\circ\sigma_1)^{-1*n} (\eta), \\
T_n & = \frac{1}{4^n} 
\sum_{f \in \calF_{n}}f^*(\eta) , 
\end{align*}
where $\calF_n = \{
\sigma_{i_1}\circ\cdots\circ\sigma_{i_n} \mid 
i_j = 1, 2, \; 1 \leq j \leq n\}$. 
Since $\sigma_1^2 = \sigma_2^2 = \id$, we have 
\[
T_{2n} = 
\frac{1}{4^{2n}} 
\sum_{i=0}^n 
\binom{2n}{i} \left\{
(\sigma_2\circ\sigma_1)^{*(n-i)}(\eta) + 
(\sigma_2\circ\sigma_1)^{-1*(n-i)}(\eta) 
\right\} 
= \frac{1}{4^{2n}} \sum_{i=0}^n \binom{2n}{i} 
\lambda^{n-i} S_{n-i}.
\]
From the definition, as $n$ goes to $\infty$, 
$S_n$ converges to $T^{+} + T^{-}$ 
and $T_{2n}$ converges to $T$. 
Thus it is sufficient to show the following claim 
to prove the proposition. 

\begin{Claim}
Let $\{s_n\}_{n=1}^{\infty}$ be a sequence with 
$s_n \in \CC$ that converges to $s_{\infty}$ as 
$n \to \infty$. We set 
$t_{2n} = \frac{1}{4^{2n}} \sum_{i=0}^n \binom{2n}{i} 
\lambda^{n-i} s_{n-i}$. 
Then $\lim_{n\to\infty} t_{2n} = s_{\infty}$. 
\end{Claim}

Note that 
$\left(
\sqrt{\lambda}+ \frac{1}{\sqrt{\lambda}}\right)^{2n} 
= \sum_{i=0}^{2n} \binom{2n}{i} \lambda^{n-i}$. 
Since $\sqrt{\lambda}= 2+\sqrt{3}$, this shows 
$\frac{1}{4^{2n}} \sum_{i=0}^{2n} \binom{2n}{i} 
\lambda^{n-i} =1$. On the other hand, for fixed $i$ 
($0 \leq i \leq n$), we have 
$\frac{1}{4^{2n}} \binom{2n}{i} \lambda^{n-i} 
= \left(\frac{\lambda}{4^{2}}\right)^{n} \binom{2n}{i} \lambda^{-i}$ goes 
to zero as $n \to \infty$.   
Moreover, we have
$\frac{1}{4^{2n}} \sum_{i=n+1}^{2n} \binom{2n}{i} 
\lambda^{n-i}
< \frac{1}{4^{2n}} \sum_{i=n+1}^{2n} \binom{2n}{i} 
< \frac{1}{4^{2n}} (1+1)^{2n}$, and thus 
$\frac{1}{4^{2n}} \sum_{i=n+1}^{2n} \binom{2n}{i} 
\lambda^{n-i}$ goes 
to zero as $n \to \infty$.   
Combining these estimates, we obtain the claim.
\QED

\subsection{Question about distribution of small points}
\label{subsec:question}
Let $(X; f_1, \cdots, f_k)$ be a
dynamical system of $k$ morphisms over a number field $K$ 
associated with a line bundle $L$ of degree $d > k$. 

In \S\ref{sec:height:subvar}, 
we have defined the canonical heights of subvarieties,  
and in \S\ref{sec:metrics:currents}
we have defined the invariant current $T$ on $X(\CC)$. 
Here, we would like to pose a question 
regarding distribution of small points. 

A sequence $\{x_n \}_{n=1}^{\infty}$ of $X(\overline{K})$
is said to be generic if for any closed 
subscheme $Y \subsetneq X$ there exists $n_0$ 
such that $x_n \not\in Y$ for $n \geq n_0$. 
A sequence $\{x_n\}_{n=1}^{\infty}$ of $X(\overline{K})$
is called a sequence of small points if 
$\lim_{n \to \infty} \widehat{h}_{L, \calF} (x_n) = 0$. 
For $x \in X(\overline{K})$, 
let $O(x)$ denote the Galois orbit of $x$ over $K$. 
Since $T$ admits a locally continuous potential 
by Theorem~\ref{thm:invariant:currents}, 
we can define a probability measure 
$\displaystyle{\mu : = \frac{1}{\chern_1(L)^{\dim X}} 
T \wedge \cdots \wedge T}$ ($\dim X$-times). 

\begin{Question}[distribution of small points]
\label{question}
Suppose there exists a generic sequence 
$\{x_n\}_{n=1}^{\infty}$ of $X(\overline{K})$
of small points. Then does 
$\frac{1}{\# O(x_n)} \sum_{y \in O(x_n)} \delta_y$ 
converge weakly to $\mu$ as $n$ goes to $\infty$? 
\end{Question}

We note that Szpiro, Ullmo and Zhang \cite{SUZ} 
proved equiditribution of small points for abelian varieties, 
Bilu \cite{Bi} for algebraic tori, 
and Chambert-Loir \cite{CL} for certain semi-abelian varieties. 

We would like to point out a similarity between the above question 
and the theorem by Briend and Duval for morphisms 
$f: \PP^N \to \PP^N$ of degree $d \geq 2$ over $\CC$.  
(Here we list only a part of their results.) 

\begin{Theorem}[Briend and Duval]
\label{thm:BD}
\begin{enumerate}
\item[(1)]\textup{(}Distribution of repelling 
periodic points of $f$.\textup{)}\,\,
The measure 
$\frac{1}{d^{kn}}
\sum_{f^n(y)=y, \,\text{$y$ is repelling}} 
\delta_y$
converge weakly to $\mu$ as $n$ goes to $\infty$.  
\item[(2)]\textup{(}Distribution of preimages of points 
outside the exceptional set.\textup{)}\,\,
The exceptional set $E$ of $f$ is defined as the biggest proper 
algebraic set of $\PP^N$ that is totally invariant by $f$. 
Then for any $a \in \PP^N(\CC)\setminus E$, 
$\frac{1}{d^{kn}}
\sum_{y \in \PP^N(\CC), f^n(y)=a} 
\delta_y$
converge weakly to $\mu$ as $n$ goes to $\infty$.  
\end{enumerate}
\end{Theorem}

Suppose $f$ is defined over a number field $K$. 
Then, for any periodic point $x$, $\widehat{h}_{\OO_{\PP^N}(1), f} (x) =0$. 
Also for $a \in \PP^N(\overline{K})\setminus E$, 
consider the set $B = \{ x \in X(\overline{K}) 
\mid f^{n}(x) = a \,\,\text{for some $n \geq0$} \}$. 
Note $\widehat{h}_{\OO_{\PP^n}(1), f} (x) = \frac{1}{d^n} 
\widehat{h}_{\OO_{\PP^N}(1), f}(a)$ if $f^{n}(x) = a$. 
Since $a \not\in E$, we can choose a generic 
sequence $\{x_n\}$ of $\PP^N(\overline{K})$ of small points 
with $x_n \in B$ for every $n \geq 0$. Thus 
Question~\ref{question} might be seen as an arithmetic analogue 
of Theorem~\ref{thm:BD} for $f: \PP^N \to \PP^N$ of degree $d \geq 2$. 

\begin{Remark}
When $f: \PP^1 \to \PP^1$ is a polynomial map, 
Baker and Hsia \cite{BH} have recently showed 
the distribution theorem for small points 
on $\PP^1(\overline{K})$. 
\end{Remark}

\subsection{Distribution of small points for Latt\`es examples}
\label{subsec:distribution:Lattes}
In this subsection, using the methods of \cite{SUZ} and \cite{CL}, 
we would like to briefly remark that Question~\ref{question} is true 
for Latt\`es examples. 

\subsubsection{Quick Review of Latt\`es examples}
\label{subsubsec:reveiw:Lattes}
Let $U(N)$ be the unitary group of size $N$, and    
$E(N) = \{(U, a) \mid U \in U(N), \; a\in \CC^N \}$ 
the complex motion group acting on $\Aff^N_{\CC}$. 
A subgroup $G$ of $E(N)$ is called a {\em complex crystallographic 
group} if $G$ is a discrete subgroup of $E(N)$ with compact quotient. 

Berteloot and Loeb \cite{BL} obtained a geometric characterization 
of morphisms of $\PP^N_{\CC}$ whose Green current is 
smooth and strictly positive on some non-empty (analytic) open subset. 

Recall that, by Proposition~\ref{prop:current:projective:space}, 
for a morphism $f: \PP^N_{\CC} \to \PP^N_{\CC}$ 
of degree $d \geq 2$, the invariant current in 
Theorem~\ref{thm:invariant:currents}
is nothing but the Green current of $f$. 

\begin{Theorem}[\cite{BL}, Th\'eor\`eme~1.1, Proposition~4.1]
\label{thm:Lattes}
Let $f: \PP^N_{\CC} \to \PP^N_{\CC}$ be a morphism over $\CC$ 
of degree $d \geq 2$, and $T$ the Green current 
of $f$. 
Assume that $T$ is smooth and strictly positive 
on some non-empty \textup{(}analytic\textup{)} open subset of $\PP^N_{\CC}$. 
Then there exists a complex crystallographic 
group $G$, and an affine transformation 
$D : \Aff^N_{\CC} \to \Aff^N_{\CC}$ whose linear part is $\sqrt{d} U$ 
with $U \in U(N)$ such that 
the ramified covering $\sigma: \Aff^N_{\CC} \to 
\PP^N_{\CC}$ satisfies $\sigma\circ f = D\circ\sigma$ and 
that $G$ acts on the fibers of $\sigma$ transitively. 
Moreover, $\sigma^*(T) 
= \frac{\sqrt{-1}}{2\pi} \sum_{i=1}^N dz_i \wedge \overline{dz_i}$, 
where $z_i$'s are the standard coordinate of $\Aff^N_{\CC}$. 
\end{Theorem}

Endomorphisms $f: \PP^N_{\CC} \to \PP^N_{\CC}$ 
that satisfy the assumption of Theorem~\ref{thm:Lattes} 
are called {\em Latt\`es examples}. 

\begin{Example}[Some Latt\`es examples]
\begin{enumerate}
\item[(1)]
Let $\tau \in \CC$ with $\Image(\tau) > 0$. 
Set $G_1 = \{ (\pm 1, m+n\tau) \mid m, n\in\ZZ\} \subset E(1)$,  
and let $D: \Aff^1 \to \Aff^1$ be the $[d]$-th map.  
Then the quotient $\Aff^1/G_1$ is isomorphic to $\PP^1$ 
and $D$ descends a morphism $f_1:  \PP^1 \to \PP^1$. 
\item[(2)](Ueda \cite{Ue})\; 
Consider the morphism $f_1 \times\cdots\times f_1: 
\PP^1 \times\cdots\times \PP^1 \to \PP^1 \times\cdots\times \PP^1$ 
($N$-times). 
The $N$-th symmetric group $S_N$ acts on 
$\PP^1 \times\cdots\times \PP^1$ by 
$\sigma \cdot (x_1, \cdots, x_n) = 
(x_{\sigma(1)}, \cdots, x_{\sigma(n)})$ for $\sigma \in S_N$.  
Observing that the quotient 
$(\PP^1 \times\cdots\times \PP^1)/S_N$ is 
isomorphic to $\PP^N_{\CC}$, 
one finds that 
$f_1 \times\cdots\times f_1$ descends 
a morphism $f_N: \PP^N_{\CC} \to  \PP^N_{\CC}$. 
\end{enumerate}
\end{Example}

We need the following proposition in the next subsection. 

\begin{Proposition}
\label{prop:Lattes:positive}
Let $f: \PP^N_{\CC} \to \PP^N_{\CC}$ be a Latt\`es example, 
and $T$ the Green current of $f$. 
\begin{enumerate}
\item[(1)]
There exists $c > 0$ such that 
$T \geq c\, \omega_{FS}$ as a current. 
\item[(2)]
Let $v: \PP^N(\CC) \to \RR$ be a $C^{\infty}$-function. 
Then there exists $\varepsilon_0 > 0$ such that, 
for all $\varepsilon$ with $0 < \varepsilon \leq \varepsilon_0$, 
$T + \varepsilon d d^c v$
is positive as a current. 
\end{enumerate}
\end{Proposition}

\Proof
Let $\sigma: \Aff^N_{\CC} \to \PP^N_{\CC}$ be the 
ramified covering as in Theorem~\ref{thm:Lattes}. 
Since $\sigma^*\omega_{FS}$ is a $C^{\infty}$ bounded 
$(1,1)$-form on $\Aff^N(\CC)$, there exists $c > 0$ such that 
$\frac{\sqrt{-1}}{2\pi} \sum_{i=1}^N dz_i \wedge \overline{dz_i}
- c\, \sigma^*\omega_{FS}$ is a strictly positive $(1,1)$-form. 
Then considering the pull-back by $\sigma$ and 
observing $\sigma^* T = \frac{\sqrt{-1}}{2\pi} 
\sum_{i=1}^N dz_i \wedge \overline{dz_i}$, 
we find $T \geq c\, \omega_{FS}$. 
The assertion (2) follows from (1). 
\QED

\subsubsection{Distribution of small points for Latt\`es examples}
\label{subsubsec:distribution:Lattes}
We follow \cite{SUZ} and \cite{CL}. In this subsection, 
we use the terminology in \cite{Zh} freely. 

\begin{Theorem}[\cite{SUZ}, Th\`eor\'eme~3.1, \cite{CL}, Proposition~6.2]
\label{thm:SUZ:CL}
Let $X$ be a projective variety defined over a number 
field $K \subset \CC$, and $L$ an ample line bundle on $X$. 
Let $\Vert\cdot\Vert$ be an ample adelic metric on $L$. 
Assume that the height of $X$ is zero with respect to 
$(L, \Vert\cdot\Vert)$. 
Let $\{x_n\}_{n=1}^{\infty}$ 
be a generic sequence of small points in $X(\overline{K})$. 
Let $v: X(\CC) \to \RR$ is a continuous function such that, 
for all sufficiently small $\varepsilon > 0$, 
$c_1(L, \Vert\cdot\Vert) + \varepsilon d d^c v$ 
is positive as a current on $X(\CC)$. 
Then we have 
\[
\liminf_{n \to \infty} \frac{1}{\# O(x_n)} \sum_{y \in O(x_n)} v(y)
\geq \int_{X} v 
\frac{\chern_1(L, \Vert\cdot\Vert)^{\dim X}}{\chern_1(L)^{\dim X}}. 
\]
\end{Theorem}

Using Theorem~\ref{thm:SUZ:CL}, we find that Question~\ref{question} 
is true for Latt\`es examples. 

\begin{Theorem}
\label{thm:distribution:Lattes}
Let $f: \PP^N \to \PP^N$ be a Latt\`es example 
defined over a number field $K \subset \CC$,
$T$ the Green current of $f$, and 
$\mu = \wedge^{N} T$. 
Let $\{x_n\}_{n=1}^{\infty}$ 
be a generic sequence of small points in $X(\overline{K})$. 
Then 
$\frac{1}{\# O(x_n)} \sum_{y \in O(x_n)} 
\delta_y$ converges weakly to $\mu$. 
\end{Theorem}

\Proof
We fix an isomorphism $\varphi: \OO_{\PP^N}(d) \simeq 
f^*(\OO_{\PP^N}(1))$, where $d$ is the degree of $f$. 
Fix a $C^{\infty}$ model $(\calX, \overline{\calL})$ 
of $(\PP^N, \OO_{\PP^N}(1))$ such that 
$\calL$ is ample and $\chern_1(\overline{\calL}_{\CC})$ 
is strictly positive on $\PP^N(\CC)$. Then, 
by \cite{Zh},~(2.3), we inductively obtain a sequence 
of $C^{\infty}$ models $(\calX_n, \overline{\calL_n})$ of 
$(\PP^N, \OO_{\PP^N}(1))$, considering the pull-backs 
by $\varphi$ (cf. Theorem~\ref{thm:adelic:sequence:from:morphisms} 
and its proof). 
By \cite{Zh},~(2.3), the models $(\calX_n, \overline{\calL_n})$ 
induces an ample adelic metric on $\OO_{\PP^N}(1)$. 

In particular, by \cite{Zh},~Theorem(2.2), 
the induced metric $\Vert\cdot\Vert_{\infty}$ 
at infinity is characterized by a continuous metric on 
$\OO_{\PP^N_{\CC}}(1)$ with 
$\Vert\cdot\Vert_{\infty}^d 
= \varphi^*f^* \Vert\cdot\Vert_{\infty}$. 
By Proposition~\ref{prop:current:projective:space}, 
$\chern_1(\OO_{\PP^N_{\CC}}(1), \Vert\cdot\Vert_{\infty})$ 
is equal to the Green current $T$.  

To apply Theorem~\ref{thm:SUZ:CL}, 
we will show that the height of $\PP^N$ is zero with respect to 
$(\OO_{\PP^N}(1), \Vert\cdot\Vert)$. 
On one hand, it follows from \cite{Zh},~Theorem(2.4) that 
$h_{\OO_{\PP^N}(1), \Vert\cdot\Vert}(\PP^N) \geq 0$. 
One the other hand, since 
$(\OO_{\PP^N}(1), \Vert\cdot\Vert)$ is ample metrized line bundle, 
if follows from \cite{Zh},~Thereom(1.10) that 
\[
\sup_{Y \subsetneq \PP^N} \inf_{x \in \PP^N \setminus Y(\overline{K})} 
h_{(\OO_{\PP^N}(1), \Vert\cdot\Vert)}(x) 
\geq 
h_{(\OO_{\PP^N}(1), \Vert\cdot\Vert)}(\PP^N),  
\]
where $Y$ runs through the algebraic subsets of $X$. 
(In the notation of Theorem~\ref{thm:canonical:height:of:subvarieties}, 
$h_{(\OO_{\PP^N}(1), \Vert\cdot\Vert)} 
= \widehat{h}_{\OO_{\PP^N}(1), f}$.)
If $x \in \PP^N(\overline{K})$ is a periodic point, 
then $h_{(\OO_{\PP^N}(1), \Vert\cdot\Vert)}(x) = 0$.  
Since the support of $T$ is $\PP^N$ (cf. \cite{Sib},~Th\'eor\`eme~1.6.5), 
the set of periodic points are Zariski dense in 
$\PP^N(\overline{K})$ (cf. Theorem~\ref{thm:BD}(1)). 
Hence the left-hand-side of the above inequality is zero, 
and so $h_{(\OO_{\PP^N}(1), \Vert\cdot\Vert)}(\PP^N) 
= 0$.

Let $v: \PP^N(\CC) \to \RR$ be a continuous function. 
We need to show 
\[
\lim_{n \to \infty} \frac{1}{\# O(x_n)} \sum_{y \in O(x_n)} v(y)
= \int_{X} v \mu.
\]
Since the $C^{\infty}$-functions are dense in 
the set of continuous functions on $\PP^N(\CC)$, 
we may assume $v$ is $C^{\infty}$. 
Since, by Proposition~\ref{prop:Lattes:positive},  
$T + \varepsilon d d^c v$ is positive for 
all sufficiently small $\varepsilon > 0$, it follows from 
Theorem~\ref{thm:SUZ:CL} that 
$\liminf_{n \to \infty} \frac{1}{\# O(x_n)} \sum_{y \in O(x_n)} v(y)
\geq \int_{X} v \mu$.
Considering $-v$ instead of $v$, we find 
$\limsup_{n \to \infty} \frac{1}{\# O(x_n)} \sum_{y \in O(x_n)} v(y)
\leq \int_{X} v \mu$. 
Thus we get the assertion.
\QED

\medskip
\section{Local canonical heights}
In this section, for normal projective varieties, 
we see that there exists canonical local heights 
for closed points 
and that the canonical height in \S1.2 decomposes into a sum of 
these canonical local heights. 
In particular, this gives yet another construction of the canonical 
heights for closed points. 
The argument goes similarly as in \cite{La},~Chap.~10 and 
\cite{CS},~\S2 which treat the case $k=1$. 

\subsection{Quick review of local heights}
In this subsection, we briefly review local heights and 
their properties. For details, we refer to Lang's book 
\cite{La},~Chap.~10. 

Let $X$ be a projective variety over a number field $K$, 
and $U$ a non-empty Zariski open set of $X$. 
In this subsection, we assume $X$ is normal. 
Let $M_K$ denote the set of absolute values of $K$, 
and $M$ denote the set of absolute values on $\overline{K}$ 
extending those of $K$. 

\begin{Definition}
\begin{enumerate}
\item[(1)]
A function $\lambda: U(\overline{K}) \times M$ is said to be {\em 
$M_K$-continuous} if, for every $v \in M$, 
$\lambda_v: U(\overline{K}) \to \RR$, $x \mapsto \lambda(x,v)$
is continuous with respect to $v$-topology. 
\item[(2)]
A function $\gamma: M_K \to \RR$ is said to be 
{\em $M_K$-constant} if $\gamma(v)=0$ for all but 
finitely many $v \in M_K$. 
For $v' \in M$ that is an extension of $v \in M_K$, 
we set $\gamma(v'):=\gamma(v)$. Then $\gamma$ is extended to a function 
$\gamma: M \to \RR$, which is also said to be $M_K$-constant.
\item[(3)]
A function $\lambda: U(\overline{K}) \times M$ is said to be {\em 
$M_K$-bounded} if there is a $M_K$-constant function $\gamma$ such that 
$\vert\alpha(x,v)\vert \leq \gamma(v)$ 
for all $(x,v) \in U(\overline{K})\times M$. 
\item[(4)]
Let $D \in \Div(X)\otimes_{\ZZ}\RR$. 
A function 
\[
\lambda_{X,D}: (X\setminus\Supp(D))(\overline{K})\times M \to 
\RR
\]
is said to be a {\em local height} (or a {\em Weil local 
height function}) associated with $D$ if it has the following property: 
There are an affine covering $\{U_i\}$ of $X$, a Cartier divisor 
$\{(U_i, s_i)\}$ representing $D$ such that 
$\alpha_i(x,v) := \lambda_{X,D}(x,v)-v\circ s_i(x)$
for $x \in (U_i\setminus\Supp(D))(\overline{K})$ and $v \in M$ 
is $M_K$-bounded and $M_K$-continuous. 
\end{enumerate}
\end{Definition}

Note that, for every $s \in K(X)^*$,
\[
\lambda_{X,\zero(s)}: (X\setminus\Supp(\zero(s)))(\overline{K})
\times M \to \RR, \quad (x, v) \mapsto v\circ s(x)
\] 
is a local height function associated with $\zero(s)$. 

We list some properties of local heights that we will use later. 

\begin{Theorem}[\cite{La}]
\label{thm:properties:of:local:heights}
\begin{enumerate}
\item[(1)]
For every $D \in \Div(X)\otimes_{\ZZ}\RR$, there exists a 
local canonical height function associated with $D$. 
\item[(2)]
Let $\lambda_{X,D}$ be a local height function associated with $D$. 
Suppose there exists a proper Zariski closed set $Z$ containing 
$\Supp(D)$ such that $\lambda_{X,D}$ is $M_K$-bounded 
on $(X\setminus\Supp(Z))(\overline{K})\times M$. 
Then $D=0$, and $\lambda_{X,D}$ extends 
uniquely to a $M_K$-bounded and $M_K$-continuous 
function on $X(\overline{K})\times M$. 
\item[(3)]
If $\lambda_{X,D}$ and $\lambda_{X,D'}$ are local height functions 
associated with $D$ and $D'$ respectively, then 
$\lambda_{X,D} + \lambda_{X,D'}$ is a local height function 
associated with $D+D'$. 
\item[(4)]
Let $g: Y\to X$ be a morphism of normal projective varieties, 
$D$ an element of $\Div(X)\otimes_{\ZZ}\RR$, and 
$\lambda_{X,D}$ a local height function associated with $D$. 
Assume that $\phi(Y)$ is not contained in $\Supp(D)$. 
Then 
\[
\lambda_D\circ(g\times\id_M): 
(Y\setminus\Supp(\phi^*D))(\overline{K})\times M \to \RR
\]
is a local height function associated with $\phi^*D$. 
\end{enumerate}
\end{Theorem}

\Proof
For (1) see \cite{La},~Chap.~10, Theorem~3.5. For (2) see 
[ibid.],~Prop~2.3, Prop~1.5, Cor~2.4. 
Note that the assumption that $X$ is normal 
is enough to show (2). For (3) and (4), see 
[ibid.],~Prop~2.1 and Prop~2.6. These proofs can apply for 
$\Div(X)\otimes_{\ZZ}\RR$. 
\QED

In the remainder of this subsection, we recall how local height 
functions are related to global height functions. 

Let $\lambda_{X,D}$ be a local height function associated 
with a divisor $D$. Then, for any point $x \in (X\setminus\Supp(D))
(\overline{K})$, the {\em associated height} is defined by 
\[
h_{\lambda_{X,D}}(x) = 
\frac{1}{[L:K]}\sum_{w \in M_L}[L_w:K_v]\lambda_{X,D}(x,w),
\]
where $L$ is an extension field of $K$ with $x\in X(L)$ and 
$w$ is an extension of $v$. The height $\lambda_{X,D}$ is defined for 
all $x \in X(\overline{K})$. Indeed, for any point 
$x\in X(\overline{K})$, there exists a rational function 
$s$ such that $x \not\in \Supp(D - \zero(s))$: Then 
we define 
\[
h_{\lambda_{X,D}}(x) := h_{\lambda_{X,D- \zero(s)}}(x). 
\]
By the product formula, this value does not depend on the choice 
of $s$.   

\begin{Theorem}[\cite{La}]
\label{thm:global:height:accociated:local:height}
Let $h_{X, \OO_X(D)}: X(\overline{K}) \to \RR$ be 
a height function corresponding to $\OO_X(D)$ 
defined in \S\ref{subsec:quick:review:of:heights}. Then 
$h_{\lambda_{X,D}} = h_{X, \OO_X(D)} + O(1)$.
\end{Theorem}

\Proof
For its proof, we refer to \cite{La},~Chap.~10, \S4.
\QED

\subsection{Construction of local canonical heights}
\label{subsec:construction:of:local:canonical:heights}
Let $X$ be a normal projective variety over a number field $K$, 
and $f_1, \cdots, f_k: X \to X$ be morphisms over $K$.  
Set $\calF=\{f_1, \cdots, f_k\}$ as before. 
Let $E$ an element of $\Div(X)\otimes_{\ZZ}\RR$. 
Assume that $f_i(X)$ is not contained in $\Supp(E)$ for each $i$, 
and that $f_1^*E + \cdots + f_k^*E \sim {d} E$ 
with $d > k$. Take $s \in \Rat(X)\otimes\RR$ with 
\[
f_1^*E + \cdots + f_k^*E = d E + \zero(s) 
\quad(\in \Rat(X)\otimes\RR). 
\]

\begin{Theorem}
\label{thm:construction:of:local:canonical:heights}
Let the notation and assumption be as above. 
Then there exists a unique function 
\[
\widehat{\lambda}_{X, E, \calF, s}: 
(X\setminus\Supp(E))(\overline{K})\times M \to \RR
\]
with the following properties\textup{:}
\begin{enumerate}
\item[(1)]
$\widehat{\lambda}_{X, E, \calF, s}$ is a Weil local height 
associated with $E$. 
\item[(2)]
For any $x \in \left(X\setminus\left(
\Supp(E)\cup\Supp(f_1^*E_1)\cup\cdots\cup\Supp(f_k^*E)
\right)\right)$ 
and any $v \in M$, 
\[
\sum_{i=1}^k 
\widehat{\lambda}_{X, E, \calF, s}(f_i(x), v)=
d \widehat{\lambda}_{X, E, \calF, s}(x, v) + v(s(x)). 
\]
\end{enumerate}
\end{Theorem}

We first prove the following lemma. 

\begin{Lemma}
\label{lemma:widehat:gamma}
Let ${\mathfrak X}$ be a topological space, 
$f_1, \cdots, f_k: {\mathfrak X} \to {\mathfrak X}$ continuous maps, and 
$\gamma: {\mathfrak X} \to \RR$ a bounded continuous function. 
Fix a real number $d$ larger than $k$. 
Then, there exists a unique bounded continuous function 
$\widehat{\gamma}: {\mathfrak X} \to \RR$ such that 
\[
\gamma(x) = \sum_{i=1}^k \widehat{\gamma}(f_i(x)) - 
d \widehat{\gamma}(x)
\]
for any $x \in {\mathfrak X}$. Moreover, 
\[ 
\Vert \widehat{\gamma} \Vert_{\sup}  
\leq \left(\frac{2d+1}{d-k} \right) \Vert\gamma\Vert_{\sup}
\]
\end{Lemma}

\Proof
We can prove the lemma as in \cite{La},~Chap.~11 Lemma~1.2 and 
\cite{CS},~Lemma~2.1. 
Namely, for a continuous function $\delta: {\mathfrak X} \to \RR$, 
define $S\delta:{\mathfrak X} \to \RR$ by 
$S\delta(x)=\frac{1}{d}
\left[\sum_{i=1}^k\delta(f_i(x))-\gamma(x)\right]$. 
\begin{Claim}
$\{S^l\gamma\}_{l=0}^{\infty}$ is a Cauchy sequence 
with respect to the sup norm.
\end{Claim}
Note that for bounded continuous functions $\delta_1, \delta_2: 
{\mathfrak X} \to \RR$, $\Vert S\delta_1 - S\delta_2\Vert_{\sup}\leq
\frac{k}{d}\Vert\delta_1 - \delta_2\Vert_{\sup}$. 
Then  $\Vert S^l\delta_1 - S^l\delta_2\Vert_{\sup}\leq
\left(\frac{k}{d}\right)^l \Vert\delta_1 - \delta_2\Vert_{\sup}$. 
Since $\sum_{l=1}^{\infty} \left(\frac{k}{d}\right)^l < \infty$, 
$\{S^l\gamma\}_{l=0}^{\infty}$ is a Cauchy sequence. 

Set $\widehat{\gamma}(x):=\lim_{l \to\infty} S^l\gamma(x)$ for 
$x \in {\mathfrak X}$. Then $\widehat{\gamma}:{\mathfrak X} \to \RR$ is a continuous 
function. It follows from $S\widehat{\gamma}(x)=\widehat{\gamma}(x)$ 
that $\gamma(x) = \sum_{i=1}^k \widehat{\gamma}(f_i(x)) - 
d \widehat{\gamma}(x)$. 
On the other hand, we have 
\begin{align*}
\Vert S^l\gamma \Vert_{\sup}
& \leq \Vert S^l\gamma - \gamma\Vert_{\sup} + \Vert\gamma\Vert_{\sup} 
\leq \sum_{\alpha=0}^{l-1}\Vert S^{\alpha+1}\gamma - 
S^{\alpha}\gamma \Vert_{\sup} +  \Vert\gamma\Vert_{\sup} \\
& \leq \sum_{\alpha=0}^{l-1} \left(\frac{k}{d}\right)^{\alpha} 
\Vert S\gamma - \gamma\Vert_{\sup} + \Vert\gamma\Vert_{\sup}
\leq \frac{d}{d-k}\Vert S\gamma-\gamma\Vert_{\sup} + \Vert\gamma\Vert_{\sup}. 
\end{align*}
Since 
\begin{align*}
\Vert S\gamma - \gamma\Vert_{\sup} 
= \sup_{x \in {\mathfrak X}} \left\vert 
\frac{1}{d}\sum_{i=1}^k \gamma(f_i(x)) - 
\left( \frac{1}{d} + 1\right) \gamma(x)\right\vert \\
\leq \left(
\frac{k}{d} + \frac{1}{d} + 1 \right) \Vert\gamma\Vert_{\sup}
= \frac{d+k+1}{d} \Vert\gamma\Vert_{\sup}, 
\end{align*}
we have 
$\Vert S^l{\gamma} \Vert_{\sup}  
\leq \left(\frac{2d+1}{d-k}\right) \Vert\gamma\Vert_{\sup}$. 
By letting $l \to \infty$, we get $\Vert \widehat{\gamma} \Vert_{\sup}  
\leq \left( \frac{2d+1}{d-k} \right) \Vert\gamma\Vert_{\sup}$. 
Finally we show the uniqueness of $\widehat\gamma$. 
Indeed suppose $\widehat\gamma_1, \widehat\gamma_2$ satisfy 
$\gamma(x) = \sum_{i=1}^k \widehat\gamma_j(f_i(x)) 
- d \widehat\gamma_j(x)$ for $j=1, 2$. 
Then $\sum_{i=1}^k (\widehat\gamma_1 - \widehat\gamma_2)(f_i(x))
= d (\widehat\gamma_1 - \widehat\gamma_2)(x)$. 
Then we have $\Vert \widehat\gamma_1 - \widehat\gamma_2\Vert_{\sup} 
\leq \frac{k}{d} \Vert \widehat\gamma_1 - \widehat\gamma_2\Vert_{\sup}$. 
Since $d > k$, we have $\widehat\gamma_1 = \widehat\gamma_2$. 
\QED

\medskip
{\sl Proof of Theorem~\ref{thm:construction:of:local:canonical:heights}.}\quad
Take a local height function $\lambda_{X,E}$. 
Since, by the assumption, $f_i(X)$ is not contained in $\Supp(E)$ 
for each $i$, $\lambda_{X, E}\circ(f_i\times \id_M)$ is a local 
height function associated with $f_i^*(E)$ by 
Theorem~\ref{thm:properties:of:local:heights}(4).  
We put 
$\lambda_{X, f_i^*E} := \lambda_{X,E} \circ (f_i\times \id_M)$. 
Set $Z:= \Supp(E)\cup\Supp(f_1^*E)\cup\cdots\cup\Supp(f_k^*E)$. 
Since $f_1^*E + \cdots + f_k^*E = d E + \zero(s)$, 
\[
\gamma(x,v) := \sum_{i=1}^k \lambda_{X, f_i^*E}(x,v)
- d\lambda_{X,E}(x,v) - v(s(x))
\]
is an $M_K$-bounded and $M_K$-continuous function on 
$(X\setminus Z)(\overline{K})\times M$. 
By Theorem~\ref{thm:properties:of:local:heights}(2), 
$\gamma$ extends to an $M_K$-bounded and $M_K$-continuous 
function on $X(\overline{K})\times M$. Then, 
by Lemma~\ref{lemma:widehat:gamma}, there exists an $M_K$-bounded and 
$M_K$-continuous function 
$\widehat{\gamma}:X(\overline{K})\times M\to\RR$ 
such that 
\[
\gamma(x,v) = \sum_{i=1}^k \widehat{\gamma}(f_i(x),v) 
- d \widehat{\gamma}(x,v) 
\]
for $x\in X(\overline{K})$ and $v\in M$. 

We define a function $\widehat{\lambda}_{X,E, \calF, s}: 
(X\setminus\Supp(E))(\overline{K})\times M \to \RR$ to be 
\[
\widehat{\lambda}_{X,E, \calF, s}(x,v) 
:= \lambda_{X,E}(x,v) - \widehat{\gamma}(x,v). 
\]
Let us see $\widehat{\lambda}_{X,E, \calF, s}$ 
enjoys the properties (1) and (2). Since $\widehat{\gamma}$ is an 
$M_K$-bounded and $M_K$-continuous function, 
$\widehat{\lambda}_{X,E, \calF, s}$ is a local height function 
associated with $E$. Moreover, for $x \in (X\setminus Z)(\overline{K})$ 
and $v \in M$, we have 
\begin{align*}
\sum_{i=1}^k & \widehat{\lambda}_{X, E, \calF, s}(f_i(x), v) 
 = \left( \sum_{i=1}^k \lambda_{X,E}(f_i(x), v) \right)
   - \left( \sum_{i=1}^k \widehat{\gamma}(f_i(x),v) \right) \\
& = \left( d \lambda_{X,E}(x,v) + v(s(x)) + \gamma(x,v) \right)
- \left( d\widehat{\gamma}(x,v) + \gamma(x,v) \right) \\
& = d\left( \lambda_{X,E}(x,v) - \widehat{\gamma}(x,v)\right)
+ v(s(x)) \\
& = d \widehat{\lambda}_{X,E, \calF, s}(x,v) +  v(s(x)).
\end{align*}

Next we see the uniqueness of $\widehat{\lambda}_{X,E, \calF, s}$. 
Suppose $\widehat{\lambda}_{X,E, \calF, s}$ and 
$\widehat{\lambda}'_{X,E, \calF, s}$ are two functions satisfying 
(1) and (2). Set 
$\delta = \widehat{\lambda}_{X,E, \calF, s} - 
\widehat{\lambda}'_{X,E, \calF, s}$. 
By Theorem~\ref{thm:properties:of:local:heights}(2), 
$\delta$ extends to an $M_K$-bounded and $M_K$-continuous function 
on $X(\overline{K})\times M$. Moreover, 
\[
\sum_{i=1}^k\delta(f_i(x),v)= d\delta(x,v)
\]
holds for $x \in (X\setminus Z)(\overline{K})$ and $v \in M$, 
and hence for all $x \in X(\overline{K})$ and $v \in M$ by 
$M_K$-continuity. Take $M_K$-constant 
$\gamma: M \to \RR$ with $\vert\delta(x,v)\vert \leq \gamma(v)$. 
Then 
\[
\vert\delta(x,v)\vert
\leq \left\vert \frac{1}{d^l} 
\sum_{f \in \calF_l}\delta(f(x),v)\right\vert
\leq \frac{k^l}{d^l}\gamma(v)
\to 0 \quad (l \to \infty),
\]
hence $\delta(x,v) = 0$ for all $x \in X(\overline{K})$ 
and $v \in M$. This shows the uniqueness.
\QED 

\subsection{Decomposition of canonical heights into local canonical heights}
In this subsection, we show that the associate height of 
$\widehat{\lambda}_{X, E, \calF, s}$ is the canonical height 
$\widehat{h}_{\OO_X(E), \calF}$ constructed in 
\S\ref{subsec:construction:of:local:canonical:heights}
coincides with the canonical height constructed in 
\S\ref{subsec:construction:of:canonical:heights}. 
In particular, this gives another construction of the canonical 
heights when $X$ is a normal projective variety. 

\begin{Theorem}
\label{thm:global:and:local:canonical:height}
Let the notation and assumption be the same as in 
Theorem~\ref{thm:construction:of:local:canonical:heights}. 
Let ${h}_{\widehat{\lambda}_{X, E, \calF, s}}$ be the 
associated height of $\widehat{\lambda}_{X, E, \calF, s}$. 
Let $\widehat{h}_{\OO_X(E), \calF}$ be the canonical height 
constructed in Theorem~\ref{thm:main}. 
Then ${h}_{\widehat{\lambda}_{X, E, \calF, s}} = 
\widehat{h}_{\OO_X(E), \calF}$. 
\end{Theorem}

\Proof
Since $\widehat{\lambda}_{X, E, \calF, s}$ 
is a local height function associated with $E$, 
${h}_{\widehat{\lambda}_{X, E, \calF, s}}
= h_{\OO_X(E)} + O(1)$ by 
Theorem~\ref{thm:global:height:accociated:local:height}.
Thus in virtue of Theorem~\ref{thm:main}, we have only to show 
\begin{equation}
\label{eqn:global:and:local:canonical:height}
{h}_{\widehat{\lambda}_{X, E, \calF, s}}(f_i(x)) 
= {d} {h}_{\widehat{\lambda}_{X, E, \calF, s}}(x)
\end{equation}
for any $x \in X(\overline{K})$. First suppose 
$x \not\in Z= \Supp(E)\cup\Supp(f_1^*E)\cup\cdots\cup\Supp(f_k^*(E))$. 
Then, the equation \eqref{eqn:global:and:local:canonical:height} 
follows from 
\[
\sum_{i=1}^k \widehat{\lambda}_{X, E, \calF, s} (f_i(x),v)
={d} \widehat{\lambda}_{X, E, \calF, s}(x,v) + v(s(x))
\]
and the product formula. 
Next, suppose $z \in Z$. 
By taking a suitable rational function $t$ and setting 
$E' := E - \zero(t)$, 
we have $x \not\in \Supp(E')\cup\Supp(f_1^*(E'))\cup\cdots\cup
\Supp(f_k^*(E'))$. 
Put $s' = s t^{{d}} \prod_{i=1}^k f_i^*(t^{-1})$ in 
$\Rat(X)\otimes \RR$. 
Then, since 
$f_1^*E'+\cdots+f_k^*E =  {d}E' + \zero(s')$, 
we similarly obtain 
${h}_{\widehat{\lambda}_{X, E', \calF, s'}}(f_i(x)) 
= {d} {h}_{\widehat{\lambda}_{X, E, \calF, s'}}(x)$. 
On the other hand, 
${h}_{\widehat{\lambda}_{X, E', \calF, s'}} = 
{h}_{\widehat{\lambda}_{X, E, \calF, s}}$ 
by the product formula. Thus the equation 
\eqref{eqn:global:and:local:canonical:height}
holds for all $x \in X(\overline{K})$. 
\QED

\section*{Appendix: Finiteness results for $\calF$-periodic points}

\renewcommand{\theTheorem}{App.~\arabic{Theorem}}
\renewcommand{\theClaim}{App.~\arabic{Theorem}.\arabic{Claim}}
\renewcommand{\theequation}{App.~\arabic{Theorem}.\arabic{Claim}}
\setcounter{Theorem}{0}

Let $K$ be a number field, and $X$ a projective variety 
defined over $K$. Let $f_1, f_2, \cdots, f_k: X \cdots\to X$ 
be rational maps. 
As in \S\ref{subsec:construction:of:canonical:heights}, 
we set $\calF_0 =\{\id\}$, $\calF_l = \{f_{i_1} \circ\cdots\circ f_{i_l} 
\mid \ 1 \leq i_1, \cdots, i_l \leq k\}$ for $l\geq 1$. 
We set $\calF :=\calF_1 (=\{f_1, \cdots, f_k \})$. 
For $f \in \bigcup_{l \geq 0}\calF_l$, let $I(f)$ denote the set 
of indeterminacy of $f$. Set $V = X \setminus 
\bigcup_{f \in \bigcup_{l \geq 0}\calF_l} I(f)$. 

We say that $x \in X(\overline{K})$ is {\em $\calF$-periodic} 
if the following conditions are satisfied: 
(1) $x\in V(\overline{K})$; 
(2) $C(x) = \{f(x) \mid f \in\bigcup_{l \geq 0}\calF_l \}$ is finite;
(3) If $C' \subseteq C(x)$ satisfies $f_i(C') \subseteq C'$ 
for every $i=1, \cdots, k$, then $C' = C(x)$.    
Note that when $k=1$ and $f_1: X \to X$ is a morphism, 
$x$ is $\calF$-periodic if and only if 
$x$ is periodic with respect to $f_1$. 

In Corollary~\ref{cor:property:of:height}, 
for dynamical systems $(X; f_1, \cdots, f_k)$ 
associated with ample line bundles, we show finiteness of the number of 
finite forward-orbits of bounded degree, 
whence finiteness of the number of 
$\calF$-periodic points of bounded degree. 
Here we would like to show finiteness of $\calF$-periodic points 
of bounded degree under a milder condition. 

\begin{Theorem}
\label{thm:app}
Let $K$ be a number field, $X$ a projective variety 
over $K$, $f_1, f_2, \cdots, f_k: X \cdots\to X$ rational maps. 
Let $S$ be a subset of $V(\overline{K})$ such that 
$f_i(S) \subseteq S$ for every $i=1, \cdots, k$. 
\textup{(}For example, we can take $S$ as $V(\overline{K})$ 
itself.\textup{)} 
Assume there exists an ample line bundle 
$L \in \Pic(X) \otimes_{\ZZ} \RR$ over $X$ 
and a positive number $\epsilon >0$ such that
\[ 
\sum_{i=1}^k h_L(f_i(x)) \geq (k+\epsilon) h_L(x) + O(1) 
\]
for all $x \in S$. 
Then for any $D \in \ZZ_{> 0}$, 
\begin{equation*}
\left\{
x \in S \ \left\vert 
\begin{gathered}
\text{$[K(x): K] \leq D$,} \\
\text{$x$ is $\calF$-periodic} 
\end{gathered}
\right. \right\}
\end{equation*}
is finite. 
\end{Theorem}

\Proof
Take a height function $h_L \in F(X)$ corresponding to $L$ and fix it. 
By the assumption, there exist a constant $C$ such that, 
for all $x \in S$,
\begin{equation}
\label{eqn:app} 
\sum_{i=1}^k h_L(f_i(x)) \geq (k+\epsilon) h_L(x) - C.
\end{equation}
Suppose $x \in S$ has the finite forward 
orbit $C(x) = \{x_1, \cdots, x_n\}$. 
For each $x_j$, consider the set $\{f_1(x_j), \cdots, f_k(x_j)\}$ and 
let $a_{ij} \in \ZZ_{\geq 0}$ denote the number such that 
$x_i$ appears $a_{ij}$ times in 
$\{f_1(x_j), \cdots, f_k(x_j)\}$  
for $i=1, \cdots, n$. 
Thus we have  
\[
h_L(f_1(x_j)) + \cdots + h_L(f_k(x_j)) 
= a_{1j} h_L(x_1) + \cdots + a_{nj} h_L(x_n). 
\] 
Since $a_{ij} \geq 0$ for all $1 \leq i,j \leq n$ and 
$a_{1j} + \cdots + a_{nj} = k$ for each $j=1, \cdots, n$,  
by Lemma~\ref{lemma:app} below, there exist $c_1, \cdots, c_n \geq 0$ 
such that $(c_1, \dots, c_n) \neq 0$ and 
$\sum_{j=1}^k a_{ij} c_j = k c_i$. 

Substituting $x = x_j$ in \eqref{eqn:app}, multiplying it by $c_j$, 
and taking the sum with respect to $j=1, \cdots, n$, we get 
\[
k \left(c_1 h_L(x_1) + \cdots + c_ n h_L(x_n) \right)
\geq (k+\epsilon) \left(c_1 h_L(x_1) + \cdots + c_ n h_L(x_n) \right) 
- C \left( \sum_{j=1}^n c_j \right). 
\] 
Thus we get $c_1 h_L(x_1) + \cdots + c_ n h_L(x_n)  
\leq \frac{C}{\epsilon} \left(\sum_{j=1}^n c_j \right)$.

Take $x_{\alpha}$ for which $h_L(x_{\alpha})$ $= 
\min\{ h_L(x_{1}), \cdots, h_L(x_{n})\}$. 
Since $c_1, \cdots, c_n \geq 0$ and $\sum_{j=1}^n c_j >0$, 
we have $h_L(x_{\alpha}) \leq \frac{C}{\epsilon}$. 
By Northcott's finiteness theorem (Theorem~\ref{thm:Northcott}), 
\[
\left\{
x \in X(\overline{K}) \ \left\vert \
\text{$[K(x): K] \leq D$ and 
$h_L(x) \leq \frac{C}{\epsilon}$} 
\right. \right\}
\]
is finite.   
By the definition of $\calF$-periodic points, 
$x \in C(x_{\alpha})$. 
Thus
\[
\left\{
x \in S \ \left\vert \
\begin{gathered}
\text{$[K(x): K] \leq D$} \\
\text{$x$ is $\calF$-periodic}
\end{gathered}
\right. \right\}
\subset 
\left\{
C(y) \ \left\vert \
\begin{gathered}
\text{$y \in S$, $[K(y): K] \leq D$,} \\
\text{and $h_L(y) \leq \frac{C}{\epsilon}$} \\
\text{$C(y)$ is finite} 
\end{gathered}
\right. \right\}
\]
is finite. 
\QED

\begin{Lemma}
Let $k, n$ be positive integers. 
\label{lemma:app}
\renewcommand{\labelenumi}{(\arabic{enumi})}
\begin{enumerate}
\item 
Let $A = (a_{ij}) \in M(n,n, \RR)$ be a real valued $(n,n)$ matrix. 
Assume that $a_{ij} \geq 0$ for all $1\leq i,j \leq n$ and 
$\sum_{i=1}^n a_{ij} = k$ for every $j =1, \cdots, n$. 
Then there exists a non-zero vector
$c={}^t (c_1, \cdots, c_n) \in \RR^n$ such that  
$c_i \geq 0$ for every $i=1, \cdots, n$ and $Ac=kc$.
\item 
Assume further that $a_{ij} > 0$  for all $1\leq i,j \leq n$. 
Then $c \in \RR^n$ as above is unique up to positive scalars.
\end{enumerate}
\end{Lemma}

\Proof
We only sketch a proof. Since (1) follows from (2) 
by taking a suitable limit, we will prove (2). 
Since $\det(A - kI_n)=0$, there exists a non-zero vector 
$c= {}^t (c_1, \cdots, c_n)$ with 
$(A - kI_n) c =0$ . Without loss of generality, 
we may assume $c_n = 1$. 
Define a $(n-1, n-1)$ matrix $B = (b_{ij})$ 
by $b_{ii} = k- a_{ii} \;(1 \leq i \leq n-1)$ 
and $b_{ij} = - a_{ij} \;(1 \leq i , j \leq n-1, i \ne j)$. 
Then by the following claim, which can be shown by the induction 
on the size of matrix, $B$ is invertible and 
every component of $B^{-1}$ is non-negative. Since 
${}^t (c_1, \cdots, c_{n-1}) = B^{-1}\;{}^t 
(a_{1n}, \cdots, a_{n-1 n})$, we have $c_i \geq 0$ 
for $i=1, \cdots, n-1$. 
\begin{Claim} 
Let $B' = (b'_{ij})$ be a 
real valued $(n-1,n-1)$ matrix. Assume that 
$b'_{ii} > 0 \;(1 \leq i \leq n-1)$, 
$b'_{ij} \leq 0 \;(1 \leq i , j \leq n-1, i \ne j)$ 
and $\sum_{i=1}^n b'_{ij} > 0 \;(1 \leq i \leq n-1)$. 
Then  $B'$ is invertible and 
every component of $B^{' -1}$ is non-negative. 
\end{Claim}
\QED

\begin{Example}[\cite{SiHenon}]
For $a \in \overline{\QQ}\setminus\{0\}$ and $b \in \overline{\QQ}$, 
a H\'enon map $\phi: \Aff^2 \to \Aff^2$ is defined 
by $\phi(x,y) = (y, y^2 + b + ax)$.  The morphism $\phi$ 
is an affine automorphism with the inverse 
$\phi^{-1}: \Aff^2 \ni (x,y) \mapsto \left(-\frac{1}{a}x^2 
- \frac{b}{a} + \frac{1}{a}y, x\right) \in \Aff^2$. 
The morphisms $\phi$ and $\phi^{-1}$    
extend to birational maps of $\PP^2$, 
which we also denote by $\phi$ and $\phi^{-1}$ respectively. 
Silverman has shown in \cite{SiHenon},~equation~(13) that 
\[
h_{\OO_{\PP^2}(1)}(\phi(x)) + h_{\OO_{\PP^2}(1)}(\phi^{-1}(x)) 
\geq  
\frac{5}{2} h_{\OO_{\PP^2}(1)}(x) + O(1)
\]
for all $x \in \Aff^2(\overline{\QQ})$. 
Note that $x \in \Aff^2(\overline{\QQ})$ is $\phi$-periodic 
if and only if $x$ is $\{\phi, \phi^{-1}\}$-periodic. 
Then, by Silverman's argument, 
or by Theorem~\ref{thm:app}, we obtain finiteness of
the number of $\phi$-periodic points in  
$\Aff^2(\overline{\QQ})$ with bounded degree 
(\cite{SiHenon},~Theorem~4.1.)
\end{Example}


\bigskip
\begin{thebibliography}{99}

\bibitem{BGS}
J.-B. Bost, H. Gillet and C. Soul\'e, 
{\em Heights of projective varieties and positive Green forms}, 
J. Amer. Math. Soc. {\bf 7} (1994), 903--1027. 

\bibitem{Bi}
Y. Bilu, 
{\em Limit distribution of small points on algebraic tori}, 
Duke Math. J. {\bf 89} (1997), 465--476.

\bibitem{BH}
M. Baker and L.-C. Hsia,  
{\em Canonical Heights, Transfinite Diameters, and Polynomial Dynamics}, 
math.NT/0305181. 

\bibitem{BL}
Berteloot\ and\ J.-J. Loeb, 
{\em Une caract\'erisation g\'eom\'etrique des exemples de 
Latt\`es de ${\Bbb P}\sp k$},
Bull. Soc. Math. France {\bf 129} (2001), 175--188.

\bibitem{Ca}
S. Cantat, 
{\em Dynamique des automorphismes des surfaces $K3$}, 
Acta Math. {\bf 187} (2001), 1--57. 

\bibitem{CL}
A. Chambert-Loir, 
{\em Points de petite hauteur sur les vari\'et\'es semi-ab\'eliennes}, 
Ann. Sci. \'Ecole Norm. Sup. {\bf 33} (2000), 789--821.

\bibitem{CS}
G. Call and H. Silverman, 
{\em Canonical heights on varieties with morphisms}, 
Compositio Math. 89 (1993), 163--205. 

\bibitem{G}
W. Gubler, 
{\em H\"ohentheorie}, With an appendix by J\"urg Kramer, 
Math. Ann. {\bf 298} (1994), 427--455.

\bibitem{HM}
A. Hinkkanen and G. J. Martin,
{\em The dynamics of semigroups of rational functions I}, 
Proc. London Math. Soc. {\bf 73} (1996), 358--384.

\bibitem{HP}
H. Hubbard and P. Papadopol, 
{\em Superattractive fixed points in $ C\sp n$}, 
Indiana Univ. Math. J. {\bf 43} (1994), 321--365.

\bibitem{La}
S. Lang, 
Fundamentals of Diophantine geometry, 
Springer, 1983.

\bibitem{Mai}
V. Maillot, 
G\'eom\'etrie d'Arakelov des vari\'et\'es 
toriques et fibr\'es en droites int\'egrables, 
M\'em. Soc. Math. Fr. No. 80, 2000. 

\bibitem{Maz}
B. Mazur, 
{\em The topology of rational points}, 
Experiment. Math. {\bf 1} (1992), 35--45.

\bibitem{Mc}
C. T. McMullen,
{\em Dynamics on $K3$ surfaces: Salem numbers and Siegel disks}, 
J. Reine Angew. Math. {\bf 545} (2002), 201--233. 

\bibitem{Mo1}
A. Moriwaki, 
{\em Arithmetic height functions over finitely generated fields}, 
Invent. Math. {\bf 140} (2000), 101--142.

\bibitem{Mo}
A. Moriwaki, 
{\em The canonical arithmetic height of subvarieties of 
an abelian variety over a finitely generated field}, 
J. reine angrew. Math. {\bf 530} (2001), 33--54.

\bibitem{Ph}
P. Philippon, 
{\em Sur des hauteurs alternatives}, 
Math. Ann. {\bf 289} (1991), 255--283.

\bibitem{Sib}
N. Sibony, 
{\em Dynamique des applications rationnelles de $\PP^k$.}
in Dynamique et g\'eom\'etrie complexes (Lyon, 1997), 
97--185, Soc. Math. France, 1999. 

\bibitem{SiHt}
J. Silverman, 
{\em The theory of height functions}, in Arithmetic geometry,
Springer (1986), 151--166.

\bibitem{SiK3}
J. Silverman, 
{\em Rational points on $K3$ surfaces: a new canonical height}, 
Invent. Math. {\bf 105} (1991), 347--373.

\bibitem{SiHenon}
J. Silverman, 
{\em Geometric and arithmetic properties of the H\'enon map},
Math. Z. {\bf 215} (1994), 237--250.

\bibitem{So}
C. Soul\'e et al, 
{\em Lectures on Arakelov Geometry},
Cambridge Univ. Press, Cambridge, 1992.

\bibitem{SUZ}
L. Szpiro, E. Ullmo\ and\ S. Zhang, 
{\em \'Equir\'epartition des petits points}, 
Invent. Math. {\bf 127} (1997), 337--347.

\bibitem{Ue}
T. Ueda, 
{\em Fatou sets in complex dynamics on projective spaces}, 
J. Math. Soc. Japan {\bf 46} (1994), 545--555.

\bibitem{Wa}
L. Wang, 
{\em Rational points and canonical heights on $K3$-surfaces 
in $ \PP\sp 1\times \PP\sp 1\times \PP\sp 1$}, 
in Recent developments in the inverse Galois problem (Seattle, WA, 1993), 
Contemp. Math. 186 (1995), 273--289.

\bibitem{Zh}
S. Zhang,
{\em Small points and adelic metrics},
J. Algebraic Geom. {\bf 4} (1995), 281--300.

\end{thebibliography}
\end{document}